
\newcount\secno
\newcount\prmno
\def\section#1{\vskip1truecm
               \global\def\currenvir{section}
               \global\advance\secno by1\global\prmno=0
               {\bf \number\secno. {#1}}
               \smallskip}

\def\subsection{\global\def\currenvir{subsection}
                \global\advance\prmno by1
               \smallskip  \ind{ (\number\secno.\number\prmno) }}
\def\subsec{\global\def\currenvir{subsection}
                \global\advance\prmno by1\smallskip
                { (\number\secno.\number\prmno)\ }}

\def\proclaim#1{\global\advance\prmno by 1
                {\bf #1 \the\secno.\the\prmno$.-$ }}

\long\def\th#1 \enonce#2\endth{%
   \medbreak\proclaim{#1}{\it #2}\global\def\currenvir{th}\smallskip}

\def\rem#1{\global\advance\prmno by 1
{\it #1} \the\secno.\the\prmno$.-$ }

\magnification 1250 \pretolerance=500 \tolerance=1000
\brokenpenalty=5000 \mathcode`A="7041 \mathcode`B="7042
\mathcode`C="7043 \mathcode`D="7044 \mathcode`E="7045
\mathcode`F="7046 \mathcode`G="7047 \mathcode`H="7048
\mathcode`I="7049 \mathcode`J="704A \mathcode`K="704B
\mathcode`L="704C \mathcode`M="704D \mathcode`N="704E
\mathcode`O="704F \mathcode`P="7050 \mathcode`Q="7051
\mathcode`R="7052 \mathcode`S="7053 \mathcode`T="7054
\mathcode`U="7055 \mathcode`V="7056 \mathcode`W="7057
\mathcode`X="7058 \mathcode`Y="7059 \mathcode`Z="705A
\def\spacedmath#1{\def\packedmath##1${\bgroup\mathsurround =0pt##1\egroup$}
\mathsurround#1
\everymath={\packedmath}\everydisplay={\mathsurround=0pt}}
 \spacedmath{2pt}

\def\iso{\vbox{\hbox to .8cm{\hfill{$\scriptstyle\sim$}\hfill}
\nointerlineskip\hbox to .8cm{{\hfill$\longrightarrow $\hfill}} }}
\def\sdir_#1^#2{\mathrel{\mathop{\kern0pt\oplus}\limits_{#1}^{#2}}}

\font\eightrm=cmr8 \font\sixrm=cmr6

\def\pc#1{\tenrm#1\sevenrm}
\def\tx{\kern-1.5pt -}
\def\cqfd{\kern 2truemm\unskip\penalty 500\vrule height 4pt depth 0pt width
4pt\medbreak} 
\def\no{n\up{o}\kern 2pt}
\def\ind{\par\hskip 1truecm\relax}

\font\pal=cmsy7

\def\sp#1{{\cal S}\kern-1pt\raise-1pt\hbox{\pal P}^{}_C(#1)}

\frenchspacing
\input xy
\xyoption{all}
\input amssym.def
\input amssym
\vsize = 25truecm \hsize = 16.1truecm \voffset = -.5truecm
\parindent=0cm
\baselineskip15pt \overfullrule=0pt

\vglue 2.5truecm \font\Bbb=msbm10

\def\mapright#1{\smash{
\mathop{\longrightarrow}\limits^{#1}}}

\def\mapdown#1{\Big\downarrow
  \rlap{$\vcenter{\hbox{$\scriptstyle#1$}}$}}

\centerline{{\bf On the residue fields of Henselian valued stable
fields\footnote{$^{\ast }$}{\rm Mathematics Subject Classification:
Primary 12E15 16K20; Secondary 12F10 12J10}}}
\bigskip

\centerline{I.D. Chipchakov\footnote {$^{\ast }$}{Partially
supported by Grants MM1106/2001 and MI-1503/2005 of the Bulgarian
Foundation for Scientific Research.}}
\par
\medskip
\centerline{Institute of Mathematics, Bulgarian Academy of Sciences}
\par
\centerline{Acad. G. Bonchev Str., bl. 8, 1113 Sofia, Bulgaria;}
\par
\centerline{e-mail: chipchak@math.bas.bg}
\par
\vskip1.truecm Abstract. Let $(K, v)$ be a Henselian valued field
satisfying the following conditions, for a given prime number $p$:
(i) central division $K$-algebras of (finite) $p$-primary dimensions
have Schur indices equal to their exponents; (ii) the value group
$v(K)$ properly includes its subgroup $pv(K)$. The paper shows that
if $\widehat K$ is the residue field of $(K, v)$ and $\widehat R$ is
an intermediate field of the maximal $p$-extension $\widehat K
(p)/\widehat K$, then the natural homomorphism Br$(\widehat K) \to $
Br$(\widehat R)$ of Brauer groups maps surjectively the
$p$-component Br$(\widehat K) _{p}$ on Br$(\widehat R) _{p}$. It
proves that Br$(\widehat K) _{p}$ is divisible, if $p
> 2$ or $\widehat K$ is a nonreal field, and that Br$(\widehat K)_{2}$ is
of order $2$ when $\widehat K$ is formally real. We also obtain that
$\widehat R$ embeds as a $\widehat K$-subalgebra in a central
division $\widehat K$-algebra $\widehat \Delta $ if and only if the
degree $[\widehat R\colon \widehat K]$ divides the index of
$\widehat \Delta $.

\vskip1.truecm \centerline{{\bf Introduction}}
\par
\medskip
This paper is devoted to the study of central division algebras and
Brauer groups of fields pointed out in the title. Let us note that a
field $E$ is said to be stable, if the Schur index ${\rm ind}(A)$ of
each finite-dimensional central simple $E$-algebra $A$ equals the
exponent exp$(A)$, i.e. the order of the similarity class $[A]$ of
$A$ in the Brauer group Br$(E)$. We say that $E$ is absolutely
stable, if its finite extensions are stable fields. Suppose that $K$
is a field with a Henselian valuation $v$ (see (2.1)). It is easily
seen that if $K$ is perfect, the value group $v(K)$ is divisible and
char$(K) = {\rm char}(\widehat K)$, where $\widehat K$ is the
residue field of $(K, v)$, then $K$ is stable if and only if
$\widehat K$ is of the same kind. This case does not make a
valuation-theoretic interest, so we focus our attention on the one
of stable $K$ and $p$-indivisible $v(K)$ (i.e. $v(K) \neq pv(K)$),
for a given prime number $p$. As it turns out, then $\widehat K$ is
a $p$-quasilocal field, i.e. it satisfies one of the following
conditions: (i) the $p$-component Br$(\widehat K) _{p}$ of
Br$(\widehat K)$ is trivial or $\widehat K$ coincides with its
maximal $p$-extension $\widehat K (p)$ in a separable closure
$\widehat K _{\rm sep}$ of $\widehat K$; (ii) every cyclic extension
of $\widehat K$ of degree $p$ embeds as a $\widehat K$-subalgebra in
each central division $\widehat K$-algebra of index $p$. We
determine the structure of Br$(\widehat K) _{p}$ and describe the
arising close relations between central division $\widehat
K$-algebras of $p$-primary dimensions and intermediate fields of
$\widehat K (p)/\widehat K$. This allows us to find when such an
intermediate field splits a given central division $\widehat
K$-algebra. Supplemented by a description of the relations between
$\widehat K$ and the quotient group $v(K)/pv(K)$ [Ch1, Theorem 2.1],
the results of the present paper enable one to characterize some of
the basic types of Henselian valued stable fields (see [Ch1, Theorem
3.1 and Sect. 4] and [Ch3, Sect. 3]). This simplifies the process of
verifying whether a given Henselian valued field is stable (see
Proposition 4.5). Note also that our research plays an essential
role in the study of the structure of Br$(K)$ carried out in
[Ch7,8]; in particular, [Ch7, Proposition 6.2] provides a
classification, up-to an isomorphism, of the abelian groups that can
be realized as reduced parts of Brauer groups of equicharacteristic
Henselian valued absolutely stable fields with totally indivisible
value groups (i.e. $p$-indivisible, for every prime $p$).
\par
\medskip
It is known that global fields and local fields are absolutely
stable (cf. [P, Sect. 17.10] and [Re, (32.19)]). The class of
stable fields is larger and of greater diversity than the subclass
of absolutely stable fields (cf. [Ch3] and [Ch1, Corollary 2.6 and
Sect. 4]). Both are singled out by the general relations between
indices and exponents (see (1.1)), and by additional restrictions
on them reflecting the specific nature of some traditionally
interesting centres. It should be noted, however, that our present
knowledge of the stability property does not bear the character of
a unified theory but is rather a collection of largely independent
results on special fields arising mainly from number theory,
commutative algebra and the theory of algebraic surfaces (cf. [Jo;
FSa; Ar, Sect. 1; MS, (16.8)]). Similarly to other topics related
to fields, simple algebras and Brauer groups (see [Am2; Pl] and [P,
Sects. 17-20]), this draws one's attention to the
valuation-theoretic approach to this area. In particular, the
interest in Henselian valued stable fields and the choice
of the topic of this paper are motivated by the fact that their
class contains Laurent formal power series fields in one variable
over local fields, and nearly all presently known stable fields $K$
with indivisible Br$(K) _{p}$, for infinitely many $p$ (cf. [Ch1,
Corollaries 4.5$\div $4.7]).
\par
\medskip
Throughout the paper, simple algebras are assumed to be associative
with a unit and finite-dimensional over their centres, Brauer groups
of fields are considered to be additively presented, Galois groups
are viewed as profinite with respect to the Krull topology, and
homomorphisms of profinite groups are supposed to be continuous. By
a $\hbox{\Bbb Z} _{p}$-extension, we mean a Galois extension with a
Galois group isomorphic to the additive group $\hbox{\Bbb Z} _{p}$
of $p$-adic integers. For any field $E$, $d(E)$ is the class of
central division $E$-algebras, ${\cal G} _{E} := {\cal G}(E _{\rm
sep}/E)$ denotes the absolute Galois group of $E$, ${\cal P}(E)$ is the set
of those prime numbers $p$ for which $E (p) \neq E$, and $_{p} {\rm
Br}(E) = \{\delta \in {\rm Br}(E)\colon \ p\delta = 0\}$. The symbol
$\pi _{E/F}$ stands for the natural homomorphism (the scalar extension
map) of Br$(E)$ into Br$(F)$, for any field extension $F/E$. When
$F/E$ is finite and separable, the corestriction mapping Br$(F) \to
{\rm Br}(E)$ is denoted by cor$_{F/E}$. Our basic notations and
terminology concerning simple algebras, Brauer groups, valuation
theory, abstract abelian torsion groups, field extensions, Galois
theory, profinite groups and Galois cohomology are standard, as
those used, for example, in [P; J; TY; F; L1; Se1] and [Ko]. The
terms ''absolutely stable" (introduced in [B]) and ''stable closed"
(used in [Ch1,2,3]) are identical in content. Note also that a
pro-$p$-group $P$ is said to be a $p$-group of Demushkin type, if it
is infinite and the homomorphism $\varphi _{a}\colon \ H ^{1} (P,
\hbox{\Bbb F} _{p}) \to H ^{2} (P, \hbox{\Bbb F} _{p})$ mapping each
$b \in H ^{1} (P, \hbox{\Bbb F} _{p})$ into the cup-product $a \cup
b$ is surjective whenever $a \in H ^{1} (P, \hbox{\Bbb F} _{p})$ and
$a \neq 0$ (where $p$ is prime and $H ^{i} (P, \hbox{\Bbb F}
_{p})\colon \ i = 1, 2$, is the $i$-th continuous cohomology group
of $P$ with coefficients in the field $\hbox{\Bbb F} _{p}$ with $p$
elements). We call a degree of $P$ the dimension of $H ^{2} (P,
\hbox{\Bbb F} _{p})$ as a vector space over $\hbox{\Bbb F} _{p}$.
Examples of such groups and more information about them are given
at the end of Sections 3, 4 and 8.
\par
\medskip
The plan of the paper is as follows: Section 1 includes
preliminaries on fields, simple algebras and Brauer groups, used in
the sequel. Section 2 contains a necessary condition for stability
of Henselian valued fields, which allows us to turn our attention
mainly to $p$-quasilocal fields $E$ such that $p \in {\cal P}(E)$. In
Section 3, we determine the structure of Br$(E) _{p}$ as an abstract
abelian group, and prove that central division $E$-algebras of
$p$-primary dimensions are cyclic and of indices equal to their
exponents. The main result of the paper is stated as Theorem 4.1. It
shows that $\pi _{E/R}$ maps Br$(E) _{p}$ surjectively on Br$(R)
_{p}$, for any intermediate field $R$ of $E (p)/E$. Theorem 4.1 also
indicates that $R$ is a splitting field of a division algebra $D \in
d(E)$ of $p$-primary dimension if and only if the degree $[R\colon
E]$ is infinite or divisible by ind$(D)$, and that $R$ embeds in $D$ as
an $E$-subalgebra if and only if $[R\colon E]$ divides ind$(D)$. In
addition, our main result implies that the class of $p$-groups of
Demushkin type of fixed degree $d \ge 0$, which are realizable as
Galois groups of maximal $p$-extensions of fields containing
primitive $p$-th roots of unity, is closed under the formation of
open subgroups. In Sections 5 and 6 we prove that the
multiplicative group $E ^{\ast }$ equals the product of the norm
groups $N(F _{1}/E)$ and $N(F _{2}/E)$, for each pair $(F _{1}, F
_{2})$ of different extensions of $E$ in $E (p)$ of degree $p$.
This result is crucial for the proof of Theorem 4.1 presented in
Section 7 (and is obtained by a method that seems to be of
independent interest). Section 8 concentrates on residue fields of
Henselian valued absolutely stable fields with totally indivisible
value groups. Our main result on this topic shows that a nonreal
and perfect field $E$ lies in the considered class if and only if
the Sylow pro-$p$-subgroups of ${\cal G} _{E}$ are of Demushkin type
whenever $p$ is a prime number for which the cohomological
$p$-dimension cd$_{p} ({\cal G} _{E})$ of ${\cal G} _{E}$ is nonzero. We refer
the reader to [Ch3, Sect. 3], for a similar but more complete
characterization of the formally real fields of this type (which
contains a description of their absolute Galois groups, up-to an
isomorphism).
\par
\medskip
A preliminary version of this paper is contained in the preprint
[Ch5], and its main result has been announced in [Ch4]. The main
results of [Ch5] (including those referred to in [Ch1,2,3]) can be
found in Sections 2, 3 and 8 of the present paper.
\par
\vskip0.75truecm \centerline{{\bf 1. Preliminaries on simple
algebras and Brauer groups} }
\par
\medskip
In this Section, we give a brief account of some fundamental results
of the classical theory of simple algebras over arbitrary fields,
which will often be used without explicit references; a more
detailed presentation of the theory can be found, for example, in
[P; Dr1] and [J]. Let $E$ be a field and $s(E)$ the class of central
simple $E$-algebras. By Wedderburn's structure theorem (cf. [P, Sect.
3.5]), each $A \in s(E)$ is isomorphic to the full matrix ring $M
_{n}(A ^{\prime })$ of order $n$ over some $A ^{\prime } \in
d(E)$; the order $n$ is uniquely determined by $A$, and so is $A
^{\prime }$, up-to an isomorphism. Algebras $A _{1}$ and $A _{2}$
in $s(E)$ are called similar (over $E)$, if the underlying division
algebras $A _{1} ^{\prime }$ and $A _{2} ^{\prime }$ are isomorphic.
This leads to the definition of Br$(E)$ as the set of similarity
classes of $s(E)$ with the group operation induced by the tensor
product in $s(E)$. It is well-known that Br$(E)$ is an abelian
torsion group; the general relations between the structure of a
division algebra $D \in d(E)$ and the similarity class $[D] \in
{\rm Br}(E)$ are described as follows (cf. [P, Sect. 14.4]):
\par
\medskip
(1.1) (i) exp$(D)$ divides ind$(D)$ and shares with it a common set
of prime divisors;
\par
(ii) $D$ decomposes into a tensor product of central division
$E$-algebras of primary dimensions; these algebras are uniquely
determined by $D$, up-to isomorphisms.
\par
\medskip
Conversely, Brauer has shown that any arrangement of positive
integers admissible by (1.1) (i) can be realized as an
index-exponent relation for some central division algebra (cf. [P,
Sect. 19.6]). The following assertions provide a useful tool for
calculating Schur indices of central simple algebras (cf. [P,
Sect. 13.4]):
\par
\medskip
(1.2) Assume that $X$ and $Y$ are finite dimensional division
algebras over an arbitrary field $E$, and at least one of them is
contained in $d(E)$. Then:
\par
(i) The $E$-algebra $X \otimes _{E} Y$ is isomorphic to $M _{k}(T)$,
for some division $E$-algebra $T$ and some divisor $k$ of the
dimensions $[X\colon E]$ and $[Y\colon E]$. In particular, if ${\rm
g.c.d.}([X\colon E], [Y\colon E]) = 1$, then $X \otimes _{E} Y$ is a
division algebra;
\par
(ii) If $Y$ is a field, then ind$(X) \ \vert $ ind$(X \otimes _{E}
Y).[Y\colon E]$, and equality holds if and only if $Y$ embeds
$E$-isomorphically into $X$; if $Y$ is a splitting field of $X$
(i.e. $[X]$ lies in the relative Brauer group Br$(Y/E)$), then
ind$(X) \ \vert \ [Y\colon E]$.
\par
\medskip
We continue with some observations that will be applied to the
study of algebraic extensions of absolutely stable and of
$p$-quasilocal fields. Let $E/E _{0}$ be an algebraic field
extension, $A$ a finite dimensional $E$-algebra, $B$ a basis of
$A$, $\Sigma $ a finite subset of $A$, $\Sigma _{1} (B)$ the set
of structural constants of $A$ determined by $B$, and $\Sigma _{2}
(B)$ the set of coordinates of the elements of $\Sigma $ with
respect to $B$. Then the extension $E _{1}$ of $E _{0}$ generated
by the union $\Sigma _{1} (B) \cup \Sigma _{2} (B)$ is finite, and
the subring $A _{1}$ of $A$ generated by $E _{1} \cup B$ is an $E
_{1}$-subalgebra of $A$ satisfying the following (cf.
[P, Sects. 9.2 and 9.4]):
\par
\medskip
(1.3) (i) The $E$-algebras $A _{1} \otimes _{E _{1}} E$ and $A$
are isomorphic;
\par
(ii) If $A/E$ is a Galois extension and $\Sigma $ contains the
roots in $A$ of the minimal polynomial over $E $ of a given
primitive element of $A/E$, then $A _{1}/E _{1}$ is a Galois
extension and the Galois groups ${\cal G}(A _{1}/E _{1})$ and ${\cal
G}(A/E)$ are canonically isomorphic;
\par
(iii) If $A \in d(E)$, then $\Sigma $ can be chosen so that $A
_{1} \in d(E _{1})$, ${\rm exp}(A _{1}) = {\rm exp}(A)$ and ${\rm
ind}(A _{1}) = {\rm ind}(A)$.
\par
\medskip
Let $E ^{\prime }/E$ be a cyclic extension of degree $m$ and
$\sigma $ a generator of ${\cal G}(E ^{\prime }/E)$. We denote by $(E
^{\prime }/E, \sigma , \beta )$ the cyclic $E$-algebra associated
with $\sigma $ and an element $\beta \in E ^{\ast }$. This algebra
is defined as a left vector space over $E ^{\prime }$ with a basis
$1, \theta , \dots , \theta ^{m-1}$, and the multiplication
satisfying the conditions $\theta ^{m} = \beta $ and $\theta \beta
^{\prime } = \sigma (\beta ^{\prime })\theta \colon \ \beta
^{\prime } \in E ^{\prime }$. The following statements
characterize these algebras and clarify their role in the
description of the relative Brauer group Br$(E ^{\prime } /E)$
(cf. [P, Sect. 15.1]):
\par
\medskip
(1.4) (i) An algebra $B ^{\prime }$ over $E$ is isomorphic to $(E
^{\prime }/E, \sigma , b ^{\prime }  )$, for some $b ^{\prime } \in
E ^{\ast }$ if and only if $B ^{\prime } \in s(E), [B ^{\prime
}\colon E] = m ^{2}$ and $E ^{\prime }$ is $E$-isomorphic to a
maximal subfield of $B ^{\prime }$;
\par
(ii) The cyclic $E$-algebras $(E ^{\prime }/E, \sigma , c)$ and $(E
^{\prime }/E, \sigma , c ^{\prime })$ are isomorphic if and only if
$c ^{\prime }c^{-1} \in N(E ^{\prime }/E)$. Moreover, the mapping of
$E ^{\ast }$ into $s(E)$ by the rule $\lambda \to (E ^{\prime }/E,
\sigma , \lambda )\colon $
\par \noindent
$\lambda \in E ^{\ast }$, induces an isomorphism of the quotient
group $E ^{\ast }/{N(E ^{\prime }/E)}$ on Br$(E ^{\prime }/E)$.
\par
\medskip
The structure of $(E ^{\prime }/E, \sigma , \beta )$ is particularly
simple when $E$ contains a primitive $m$-th root of unity
$\varepsilon $. Then $E ^{\prime }/E$ is a Kummer extension (cf.
[L1, Ch. VIII, Theorem 10]) and there are elements $\alpha \in E
^{\ast }$ and $\xi \in E ^{\prime }$, such that $E ^{\prime }= E(\xi
)$ and $(E ^{\prime }/E, \sigma , \beta ) = E \langle \xi , \theta
\colon \ \xi ^{m} = \alpha , \theta ^{m} = \beta , \theta \xi =
\varepsilon \xi \theta \rangle $. In this case, $(E ^{\prime }/E,
\sigma , \beta )$ is called a symbol algebra and usually is denoted
by $A _{\varepsilon }(\alpha , \beta ; E)$. The general properties
of symbol algebras and their analogues of dimension $p ^{2}$ over
fields of characteristic $p > 0$, see [Se3, Ch. XIV, Sect. 1],
enable one to prove the following lemma (and Lemma 7.2) by a
separate discussion of the special cases of $p \neq {\rm char}(E)$
and $p = {\rm char}(E)$ (see, e.g. [Ch4]). For convenience of the
reader, we present here unified proofs suggested by the referee.
\par
\medskip
{\bf Lemma 1.1.} {\it Let $E$ be a field and $L$ an extension of
$E$ presentable as a compositum $L = F _{1}F _{2}$ of distinct
cyclic extensions $F _{1}$ and $F _{2}$ of $E$ of prime degree
$p$. Assume also that $F _{3}$ is an intermediate field of $L/E$,
such that $[F _{3}\colon E] = p$ and $F _{3} \neq F _{j}\colon \ j
= 1, 2$. Then $N(F _{3}/E)$ includes the intersection $N(F _{1}/E)
\cap N(F _{2}/E)$.}
\par
\medskip
Before proving the lemma, let us recall that the character group
$C(G)$ of any profinite group $G$ is the abelian group of continuous
homomorphisms of $G$ into the discrete quotient group $\hbox{\Bbb
Q}/\hbox{\Bbb Z}$ of the additive group $\hbox{\Bbb Q}$ of rational
numbers by the subgroup $\hbox{\Bbb Z}$ of integers. This is
equivalent to the standard definition of a character group used in
topological group theory, in spite of the fact that $\hbox{\Bbb
Q}/\hbox{\Bbb Z}$ is not a discrete subset of the compact group
$\hbox{\Bbb R}/\hbox{\Bbb Z}$, where $\hbox{\Bbb R}$ is the additive
locally compact group of real numbers (see [K, Ch. 7, Corollary
5.3]). Note that $C(G)$ is a torsion group, since each character of
$G$ maps it into a compact, hence a finite, subgroup of $\hbox{\Bbb
Q}/\hbox{\Bbb Z}$. Regarding $\hbox{\Bbb Q}/\hbox{\Bbb Z}$ as a
trivial $G$-module, we also identify $C(G)$ with the continuous
cohomology group $H ^{1} (G, \hbox{\Bbb Q}/\hbox{\Bbb Z})$. This
allows us to identify, for each prime $p$, the set $_{p} C(G) =
\{\chi \in C(G)\colon \ p\chi = 0\}$ with the continuous
homomorphism group ${\rm Hom} (G, \hbox{\Bbb F} _{p})$, where
$\hbox{\Bbb F} _{p}$ denotes the field with $p$ elements (and  is
viewed as a discrete additive abelian group). When $G = {\cal G}
_{E}$, for a given field $E$, we put $X _{p} (E) = \ _{p} C({\cal G}
_{E})$.
\par
\medskip
{\it Proof of Lemma 1.1.} It is clear from Galois theory that
$L/E$ is abelian and $F _{3}/E$ is cyclic. Let $X _{p} (E) = {\rm
Hom}({\cal G} _{E}, \hbox{\Bbb F} _{p})$, and for each $\chi \in X _{p}
(E)$, denote by $L _{\chi }$ the extension of $E$ in $E _{\rm
sep}$ fixed by the kernel ${\rm Ker}(\chi )$. Also, let $\chi
_{i}$ be characters in $X _{p} (E)$ such that $F _{i}$ is the
fixed field of ${\rm Ker}(\chi _{i})$, $i = 1, 2, 3$. Because $F
_{3} \subset F _{1}F _{2}$, $\chi _{3}$ lies in the span of $\chi
_{1}$ and $\chi _{2}$ in $X _{p} (E)$. Take any $c \in N(F _{1}/E)
\cap N(F _{2}/E)$ and consider the pairing $s\colon \ X _{p} (E)
\times E ^{\ast } \to $
$_{p} {\rm Br}(E)$ defined by the rule
$s(\chi , b) = (L _{\chi }/E, \sigma , b)$, where $\sigma $ maps
to $1$ in the map ${\cal G}(L _{\chi }/E) \to \hbox{\Bbb F} _{p}$ induced
by $\chi $. (If $\chi = 0$, then $L _{\chi } = E$ and $\sigma =
{\rm id}$.) By (1.4) (ii), $s(\chi _{1}, c) = s(\chi _{2},
c) = 0$. Since $s$ is $\hbox{\Bbb Z}$-bilinear (see [Se3, Ch. XIV,
Sect. 1]) and $\chi _{3}$ is an $\hbox{\Bbb F} _{p}$-linear
combination of $\chi _{1}$ and $\chi _{2}$, it follows that
$s(\chi _{3}, c) = 0$; hence, by (1.4) (ii), $c \in N(F _{3}/E)$.
\par
\medskip
When $L/E$ is an arbitrary cyclic field extension, the image
Im$(\pi _{E/L})$ of $\pi _{E/L}$ is characterized by
Teichm$\ddot u$ller's theorem (cf. [Dr1, Ch. 9, Theorem 4]) as
follows:
\par
\medskip
(1.5) For an algebra $A \in s(L)$, the following assertions are
equivalent:
\par
(i) $[A]$ lies in Im$(\pi _{E/L})$, i.e. $A$ is similar over $L$ to
$A _{0} \otimes _{E} L$, for some $A _{0} \in s(E)$;
\par
(ii) $[A]$ is fixed by the natural action of ${\cal G}(L/E)$ on
Br$(L)$, i.e. every $E$-automorphism $\psi $ of the field $L$ is
extendable to an automorphism $\bar \psi $ of $A$ (viewed as an
algebra over $E$).
\par
\medskip
The next lemma is used in Section 4 for proving that Br$(L) _{p}
\subseteq {\rm Im}(\pi _{E/L})$ in case $E$ is $p$-quasilocal
and $L/E$ is a cyclic $p$-extension.
\par
\medskip
{\bf Lemma 1.2.} {\it Let $p$ be a prime number and $H = \langle
h\rangle $ a cyclic group of order $p ^{t}$, for some $t \in
\hbox{\Bbb N}$. Then $h - 1$ is a nilpotent element of the group
ring $(\hbox{\Bbb Z}/p ^{s}\hbox{\Bbb Z}) [H]$, for any $s \in
\hbox{\Bbb N}$.}
\par
\medskip
{\it Proof.} The binomial expansion shows that $(h - 1) ^{p ^{t}} =
py$, for some $y \in (\hbox{\Bbb Z}/p ^{s}\hbox{\Bbb Z}) [H]$. Hence,
$(h - 1) ^{p ^{t}.s} = (py) ^{s} = p ^{s}y ^{s} = 0$.
\par
\medskip
The following results enable us to take in Section 6 the
crucial technical step towards proving the main result of this
paper. They are well-known consequences of Amitsur's theorem [Am1]
(see also [Roq1], [Roq2, Sect. 1] and [Am3, pp 1-3]) about the
function fields of Brauer-Severi varieties:
\par
\medskip
(1.6) (i) Every subgroup $U$ of Br$(E)$ is equal to Br$(\Lambda
_{U}/E)$, for some compositum $\Lambda _{U}$ of function fields of
Brauer-Severi varieties defined over $E$; also, $E$ is
algebraically closed in $\Lambda _{U}$ (cf. [FS, Theorem 1]);
\par
(ii) There exists a set $\{\Lambda _{n}\colon \ n \in \hbox{\Bbb
N}\}$ of extensions of $E$, such that $\Lambda _{1} = E$, and for
each index $n$, $\Lambda _{n} \subseteq \Lambda _{n+1}$, Br$(\Lambda
_{n+1}/\Lambda _{n}) = $ Br$(\Lambda
_{n})$ and $\Lambda _{n}$ is algebraically closed in $\Lambda
_{n+1}$; in particular, the union $\Lambda ^{\prime } = \cup _{n=1}
^{\infty } \Lambda _{n}$ is a field with Br$(\Lambda ^{\prime }) =
\{0\}$, in which $E$ is algebraically closed.
\par
\medskip
We end this Section by defining most of the special types of fields
used in the sequel. A field $F$ is said to be formally
real, if $-1$ is not presentable as a finite sum of elements of the
set $F ^{\ast 2} = \{f ^{2}\colon \ f \in F ^{\ast }\}$; we say
that $F$ is nonreal, otherwise. The field $F$ is called
Pythagorean, if it is formally real and $F ^{\ast 2}$ is
closed under the addition in $F$. This property is characterized as
follows (cf. [Wh, Theorem 2]):
\par
\medskip
(1.7) $F$ is Pythagorean if and only if $2 \in {\cal P}(F)$ and $F$ does
not admit a cyclic extension of degree $4$.
\par
\medskip
A field $E$ is said to be almost perfect, if every finite
extension of $E$ has a primitive element. It follows from the
general theory of algebraic extensions that if char$(E) = q > 0$,
then $[E _{1}\colon E _{1} ^{q}] = [E\colon E ^{q}]$ , for each
finite extension $E _{1}$ of $E$, where $E ^{q} _{1} = \{\alpha
_{1} ^{q}\colon \ \alpha _{1} \in E _{1}\}$ (cf. [L1, Ch. VII,
Sect. 7, Corollary 4; Ch. VIII, Sect. 9, Corollary 1]). This
implies that
\par
\medskip
(1.8) $E$ is almost perfect if and only if char$(E) = 0$ or char$(E)
= q > 0$ and $[E\colon E ^{q}]$ equals $1$ or $q$. The classes of
perfect fields and of almost perfect fields are closed under the
formation of algebraic extensions.
\par
\medskip
It is known that complete discrete valued fields with perfect
residue fields are almost perfect (cf. [E, (5.7)$\div $(5.10)]).
We say that a field $F$ is quasilocal, if every finite extension
$F _{1}/F$ satisfies the following condition:
\par
Every cyclic extension $F _{1} ^{\prime }$ of $F _{1}$ is
embeddable as an $F _{1}$-subalgebra in each $D _{1} \in d(F
_{1})$ of index divisible by $[F _{1} ^{\prime }\colon F _{1}]$.
\par
\medskip
We prove in Section 8 that $F$ is quasilocal if and only if its
finite extensions are $p$-quasilocal fields, for every prime $p$.
In particular, this is the case, if $F$ is a formally real closed
or a local field (see [L1, Ch. XI, Theorem 1] and [Se3, Ch. XIII,
Sect. 3]). Other examples of quasilocal fields can be found, for
instance, in [Se3; Ch2,3,7] and Section 8.
\par
\vskip0.75truecm {{\bf 2. A necessary condition for stability of
Henselian valued fields}}
\par
\medskip
Let $K$ be a field with a nontrivial Krull valuation $v$, $O _{K}$
the valuation ring, $\widehat K$ the residue field and $v(K)$ the
value group of $(K, v)$. We say that $v$ is Henselian, if any of
the following three conditions holds (see [R; Er] or [W], for a
proof of their equivalence):
\par
\medskip
(2.1) (i) For every monic polynomial $f \in O_{K}[X]$ whose image
$\hat f \in \widehat K [X]$ (under the natural ring homomorphism
$O _{K} [X] \to \widehat K [X]$, mapping $O _{K}$ canonically on
$\widehat K$) has a simple root $\tilde \alpha \in \widehat
K$, there is a root $\alpha \in O _{K}$ of $f$ with $\hat \alpha =
\tilde \alpha $;
\par
(ii) $v$ can be extended to a uniquely determined (up-to an
equivalence) valuation $v _{K'}$ on each algebraic extension $K
^{\prime }$ of $K$;
\par
(iii) $v$ is uniquely extendable to a valuation $v _{D}$ on each
division $K$-algebra $D$ of finite dimension $[D\colon K]$.
\par
\medskip
It is well-known that $v$ is Henselian in the following two
special cases: (i) $v$ is real-valued and $K$ is complete with
respect to the topology induced by $v$; (ii) $K$ is an iterated
Laurent formal power series field in $n \ge 1$ indeterminates and
$v$ is the standard $\hbox{\Bbb Z} ^{n}$-valued valuation of $K$.
Note also that the fulfillment of conditions (2.1) guarantees that
they are satisfied by the prolongation of $v$ on any algebraic
extension of $K$.
\par
\medskip
Assume that $v$ is Henselian and, with notation being as in (2.1)
(iii), let $\widehat D$ and $v(D)$ be the residue division ring
and the value group of $(D, v _{D})$, respectively. It is known that
$\widehat D$ is a $\widehat K$-algebra such that $[\widehat D\colon
\widehat K] \le [D\colon K]$, and $v(D)$ is a totally ordered
abelian group including $v(K)$ as a subgroup of index $e(D/K) \le
[D\colon K]$. By the Ostrowski-Draxl theorem [Dr2], $[D\colon K]$,
$[\widehat D\colon \widehat K]$ and $e(D/K)$ are related as follows:
\par
\medskip
(2.2) $[D\colon K] = [\widehat D\colon \widehat K]e(D/K)d(D/K)$,
for some integer $d(D/K)$ (called a defect of $D$ over $K$); if
$d(D/K) \neq 1$, then char$(\widehat K) = q > 0$ and $d(D/K)$ is a
power of $q$.
\par
\medskip \noindent
This, combined with (1.1) and [TY, Theorem 4.1], leads to a
complete description of the relations between Schur indices and
defects of central division algebras over Henselian valued fields.
The division $K$-algebra $D$ is said to be defectless (with respect to
$v$) if $d(D/K) = 1$, and it is called inertial if $[D\colon
K] = [\widehat D\colon \widehat K]$ and the centre $Z(\widehat D)$
of $\widehat D$ is a separable extension of $\widehat K$. By
Theorem 2.8 (a) of [JW], for every finite dimensional division
$\widehat K$-algebra $\widetilde S$ with $Z(\widetilde S)$
separable over $\widehat K$, there exists an inertial division
$K$-algebra $S$ such that $\widehat S$ is $\widehat K$-isomorphic
to $\widetilde S$. This algebra is uniquely determined by
$\widetilde S$, up- to a $K$-isomorphism (and is called an inertial
lift of $\widetilde S$ over $K$). We refer the reader to [JW], for
a systematic presentation of inertial, totally ramified, nicely
semiramified and other types of defectless central division
$K$-algebras.
\par
\medskip
The starting point for our further considerations is the following
necessary condition for stability of Henselian valued fields; its
proof has been suggested by the referee and is considerably
shorter than the one in the first version of the paper.
\par
\medskip
{\bf Proposition 2.1.} {\it Let $K$ be a stable field with a
Henselian valuation v. Then the residue field $\widehat K$ of $(K,
v)$ is also stable. Moreover, if $v(K)$ is $p$-indivisible, for
some prime number $p$, and if $\widetilde S \in d(\widehat K)$
is an algebra of $p$-primary dimension, then every cyclic extension
of $\widehat K$ of degree dividing ${\rm ind}(\widetilde S)$ is
embeddable in $\widetilde S$ as a $\widehat K$-subalgebra.}
\par
\medskip
{\it Proof}. Let $i(K) = \{\Delta \in d(K)\colon \ \Delta $
is inertial over $K\}$. By [JW, Theorem 2.8 (b)], the set IBr$(K) =
\{[\Delta ]\colon \ \Delta \in i(K)\}$ is a subgroup of Br$(K)$
and the mapping $\pi \colon \ i(K) \to d(\widehat K )$ by the
formula $\pi (S) = \widehat S\colon \ S \in i(K)$, induces a group
isomorphism $\tilde \pi \colon \ {\rm IBr}(K) \cong {\rm
Br}(\widehat K)$. In particular, $\pi $ preserves indices
and exponents, so $\widehat K$ inherits the stability of $K$.
\par
Suppose now that $pv(K) \neq v(K)$, for some prime $p$. Fix an
algebra $\widetilde S \in d(\widehat K )$ of $p$-primary
index as well as a cyclic extension $\widetilde L$ of $\widehat K$
of degree $n$ dividing ind$(\widetilde S)$, and denote by $S$ and
$L$ the inertial lifts over $K$ of $\widetilde S$ and $\widetilde
L$, respectively. So, $S \in d(K)$, and ind$(\widetilde S) = {\rm
ind}(S) = {\rm exp}(S) = {\rm exp}(\widetilde S)$, as K is stable.
Note that $L/K$ is cyclic (see [JW, page 135]), and for $b \in K
^{\ast }$ with $v(b) \not\in pv(K)$, let $V = (L/K, \sigma ,b)$,
where $\sigma$ is any generator of ${\cal G}(L/K)$. The choice of $b$
guarantees that the image of $v(b)$ in $v(K)/nv(K)$ has order $n$,
so [JW, Exercise 4.3] shows that $V$ is a nicely semiramified
division $K$-algebra with $v(V)/v(K)$ cyclic of order $n$ and
$\widehat V = \widetilde L$. Since $S$ is inertial, [JW, Theorem
5.15 (a)] and the noted property of $V$ yield that for the
underlying division algebra $D$ of $S \otimes _{K} V$, we have
exp$(D) = {\rm l.c.m.}({\rm exp}(\widetilde S), {\rm
exp}(v(V)/v(K))) = {\rm exp}(\widetilde S)$ and ind$(D) = {\rm
ind}(\widetilde S \otimes _{\widehat K} \widetilde L).n$. As the
stability of $K$ requires that ind$(D) = {\rm exp}(D)$, one
obtains from these results that ind($\widetilde S) = {\rm
exp}(\widetilde S) = {\rm ind}(\widetilde S \otimes _{\widehat K}
\widetilde L).n$. Hence, by (1.2) (ii), $\widetilde L$ embeds in
$\widetilde S$, as desired.
\par
\medskip
{\bf Corollary 2.2.} {\it Under the hypotheses of Proposition 2.1,
if $v(K)$ is totally indivisible, then every cyclic extension
$\widetilde L$ of $\widehat K$ embeds $\widehat K$-isomorphically
in each algebra $\widetilde D \in d(\widehat K)$ of index
divisible by $[\widetilde L\colon \widehat K]$.}
\par
\medskip
{\it Proof.} This can be deduced from (1.1), (1.2) and
Proposition 2.1, since $\widetilde L$ is presentable as a tensor
product over $\widehat K$ of cyclic extensions of $\widehat K$
of primary degrees (see [P, Sect. 15.3]).
\par
\medskip
{\bf Proposition 2.3.} {\it Quasilocal fields are absolutely
stable. Residue fields of Henselian valued absolutely stable
fields with totally indivisible value groups are quasilocal and
almost perfect.}
\par
\medskip
{\it Proof.} Note first that if $(K, v)$ is a Henselian valued field
and $\widetilde L/\widehat K$ is a finite extension, then there
exists an extension $L/K$ such that $[L\colon K] = [\widetilde
L\colon \widehat K]$, $v _{L} (L) = v(K)$ and the residue field of
$(L, v _{L})$ is $\widehat K$-isomorphic to $\widetilde L$ (see
(2.1), (2.2) and [L1, Ch. VII, Sects. 4 and 7]). Therefore, our second
assertion follows from Corollary 2.2 and [Ch1, Corollary 2.6
(iii)]. Let now $E$ be a quasilocal field and let $\Delta \in d(E)$ be
of prime exponent $p$. Since finite extensions of $E$ are
quasilocal, it suffices for the proof of the absolute stability of
$E$ to show that ind$(\Delta ) = p$ (see (1.1) and [A1, Ch. XI,
Theorem 3]). Suppose first that $E$ contains a primitive $p$-th
root of unity or char$(E) = p$. By the Merkurjev-Suslin theorem
and Albert's theory of $p$-algebras (cf. [MS, (16.1)] and [A1, Ch.
VII, Theorem 28]), then $\Delta $ is similar to a tensor product
of cyclic division $E$-algebras of index $p$. Therefore, $p \in
{\cal P}(E)$ and $\Delta /E$ is split by any cyclic extension of $E$ of
degree $p$, which enables one to deduce the equality ind$(\Delta )
= p$ from (1.1) and (1.2). Assume now that $p \neq {\rm char}(E)$,
$E ^{\prime }$ is an extension of $E$ obtained by adjunction of a
primitive $p$-th root of unity, and $\Delta ^{\prime } = \Delta
\otimes _{E} E ^{\prime }$. It is known that $[E ^{\prime }\colon E]$
divides $p - 1$ (cf. [L1, Ch. VIII, Sect. 3]); hence, by (1.2),
$\Delta ^{\prime } \in d(E ^{\prime })$, ind$(\Delta ^{\prime }) =
{\rm ind}(\Delta )$ and exp$(\Delta ^{\prime }) = p$. As $E
^{\prime }$ is quasilocal, this yields ind$(\Delta ) = p$.
\par
\medskip
It would be of interest to know whether a Henselian discrete
valued field $(K, v)$ is absolutely stable when $\widehat K$ is
quasilocal and almost perfect (see [Ch1, Corollary 4.6]).
\vskip0.75truecm
\centerline{\bf 3. Central division algebras of $p$-primary
dimensions and the $p$-component}
\par
\centerline{\bf of the Brauer group in the case of a
$p$-indivisible value group}
\par
\medskip
The main result of this Section sheds light on the nature of the
stability property of residue fields of Henselian valued stable
fields with $p$-indivisible value groups, and on the structure of
the $p$-components of their Brauer groups.
\par
\medskip
{\bf Theorem 3.1.} {\it Let $E$ be a $p$-quasilocal field such that
Br$(E) _{p} \neq \{0\}$, for some $p \in {\cal P}(E)$. Assume also
that $R$ is a cyclic extension of $E$ in $E (p)$, and $D \in d(E)$
is an algebra of $p$-primary index. Then $E$, $R$ and $D$ have the
following properties:}
\par
(i) $D$ {\it is a cyclic $E$-algebra and {\rm ind}$(D) = {\rm
exp}(D)$;}
\par
(ii) Br$(E) _{p}$ {\it is a divisible group unless $p = 2$ and $E$
is formally real;}
\par
(iii) {\it If $[R\colon E] = {\rm ind}(D)$, then $R$ embeds in $D$
as an $E$-subalgebra;}
\par
(iv) {\it If $[R\colon E] \le {\rm ind}(D)$, then {\rm exp}$(D
\otimes _{E} R) = {\rm exp}(D)/[R\colon E]$.}
\par
\medskip
{\it Proof.} Suppose first that $[R\colon E] = {\rm exp}(D) = p
^{n}$, for some $n \in \hbox{\Bbb N}$, fix a generator $\sigma $
of ${\cal G}(R/E)$, and denote by $R ^{\prime }$ the (unique) extension
of $E$ in $R$ of degree $p ^{n-1}$. We show that ind$(D) = p ^{n}$
and $R$ embeds in $D$ as an $E$-subalgebra. In view of (1.1) and
(1.2) (ii), this amounts to proving that $D$ is split by $R$, i.e.
$[D] \in {\rm Br}(R/E)$. Note first that the underlying division
$E$-algebra $D ^{\prime }$ of the $p$-th tensor power of $D$ over
$E$ has exponent $p ^{n-1}$, and $R ^{\prime }/E$ is cyclic with
${\cal G}(R ^{\prime }/E)$ generated by the automorphism $\sigma
^{\prime }$ of $R ^{\prime }$ induced by $\sigma $. Also, if $n =
1$, then $D$ is similar to a tensor product of algebras in $d(E)$
of index $p$ (cf. [M, Sect. 4, Theorem 2]); hence, by (1.2) and the
$p$-quasilocal property of $E$, $[D] \in {\rm Br}(R/E)$.
Proceeding by induction on $n$, one may assume further that $n \ge
2$ and $[D ^{\prime }] \in {\rm Br}(R ^{\prime }/E)$. Now, by
(1.2) and (1.4), there is an $E$-isomorphism $D ^{\prime }
\cong (R ^{\prime }/E, \sigma ^{\prime }, \lambda )$ for some
$\lambda \in E ^{\ast }$, and by [P, Sect. 15.1, Corollary b], $D
^{\prime }$ is similar to the $p$-th tensor power of the
$E$-algebra $(R/E, \sigma , \lambda )$. In view of (1.4) (ii),
this means that $D \otimes _{E} (R/E, \sigma , \lambda ^{-1})$ is
of exponent $p$ (in $s(E)$); hence, $[D \otimes _{E} (R/E, \sigma
, \lambda ^{-1})]$ and $[D]$ lie in Br$(R/E)$, which proves our
assertion.
\par
Suppose now that $[R\colon E] = p ^{n} \le {\rm exp}(D)$, put $\nu
= {\rm exp}(D)/p ^{n}$ and denote by $D _{0}$ the underlying
division $E$-algebra of the $\nu $-th tensor power of $D$ (over
$E$). It is easily verified that exp$(D _{0}) = p ^{n}$ and by the
first part of our proof, $R$ splits $D _{0}$. This, combined with
(1.1) (i) and (1.2) (ii), proves Theorem 3.1 (iv) and shows that
Br$(R/E) = \{b \in {\rm Br}(E)\colon \ p ^{n}b = 0\}$. More
precisely, arguing as above, one obtains from [P, Sect. 15.1,
Corollary b] and (1.4) (ii) that Br$(R/E)$ coincides with the
union of its cyclic subgroups of order $p ^{n}$. Now the remaining
assertions of Theorem 3.1 can be deduced from the following
statements:
\par
\medskip
(3.1) Let $E$ be a field and $p \in {\cal P}(E)$. Then:
\par
(i) For each $\Delta \in d(E)$ of $p$-primary index, there exists a
cyclic extension of $E$ in $E (p)$ of degree equal to exp$(\Delta
)$;
\par
(ii) If $p > 2$ or $E$ is nonreal, then $E (p)$ contains as a
subfield a cyclic extension of $E$ of degree $p ^{m}$, for each $m
\in \hbox{\Bbb N}$.
\par
\medskip
Statements (3.1) are obtained as consequences of the following two
lemmas.
\par
\medskip
{\bf Lemma 3.2.} {\it Let $E$ be a field and $p \in {\cal P}(E)$. Then $E
(p)$ possesses a subfield that is a $\hbox{\Bbb Z} _{p}$-extension
of $E$ unless $p = 2$ and $E$ is Pythagorean.}
\par
\medskip
For a proof of Lemma 3.2, we refer the reader to [Wh, Theorem
2]. Our next lemma is also known but its proof is included here
because of its brevity and simplicity.
\par
\medskip
{\bf Lemma 3.3.} {\it Let $E$ be a Pythagorean field. Then {\rm
Br}$(E) _{2}$ is of exponent $2$.}
\par
\medskip
{\it Proof.} As $E$ is formally real, the equation $X ^{2} + Y
^{2} = -1$ has no solution in $E$, which means that $-1 \not\in
N(E(\sqrt{-1})/E)$. Therefore, by (1.4) (ii), $A _{-1}(-1, -1; E)
\in d(E)$, whence Br$(E)_{2} \neq \{0\}$. Since $E$ is Pythagorean,
and by [MS, (16.6)], central division algebras of exponent $4$ are
similar to tensor products of cyclic algebras, it follows from (1.4)
(ii) and (1.7) that Br$(E)_{2}$ does not contain elements of order
$4$. Thus Lemma 3.3 is proved.
\par
\medskip
{\bf Remark 3.4.} (i) Lemma 3.3 is a part of known characterizations
of fields $E$ with Br$(E) _{2}$ of exponent $2$ (see, for instance,
[Ef, Theorem 3.1]). The availability of this property implies that
$E$ is formally real (not necessarily Pythagorean, e.g. the
$\hbox{\Bbb Z} _{2}$-extension of the field $\hbox{\Bbb Q}$ of
rational numbers in $\hbox{\Bbb Q} (2)$) and its totally positive
elements are presentable as sums of two squares over $E$.
\par
(ii) Note that if $E$ is a field with ${\cal G}(E(p)/E)$ a
pro-$p$-group of rank $1$ and order $\ge 3$, for some $p \in {\cal
P}(E)$, then $E (p)/E$ is a $\hbox{\Bbb Z} _{p}$-extension. Indeed,
by Galois theory, ${\cal G}(E(p)/E)$ has a unique open subgroup of index $p$.
Therefore, finite extensions of $E$ in $E (p)$ are cyclic, so our
assertion reduces to a consequence of Lemma 3.2 and statement (1.7).
\par
\medskip
Theorem 3.1 (ii) is supplemented by the following lemma.
\par
\medskip
{\bf Lemma 3.5.}  {\it A formally real field $E$ is $2$-quasilocal
if and only if $[E(2)\colon E] = 2$; when this is the case, $E$ is
Pythagorean and Br$(E) _{2}$ is of order $2$.}
\par
\medskip
{\it Proof.} Evidently, if $[E (2)\colon E] = 2$, then $(\Delta
\otimes _{E} \Delta ^{\prime }) \not\in d(E)$, for any $\Delta $ and
$\Delta ^{\prime } \in d(E)$ with ind$(\Delta ) = {\rm ind}(\Delta
^{\prime }) = 2$; also, $E$ is Pythagorean, by (1.7). This, combined
with [MS, (16.1)], (1.4) (ii) and Lemma 3.3, implies that Br$(E) _{2}
= \langle [A _{-1} (-1, -1; E)]\rangle \neq \{0\}$. Thus the latter
part and the sufficiency in the former part of the lemma are proved.
Suppose now that $E$ is $2$-quasilocal and put $B _{1} =
B(\sqrt{-1})$ (where $\sqrt{-1} \in E (2)$), for each formally real
extension $B$ of $E$ in $E (2)$. As in the proof of Lemma 3.3, one
obtains that  $[B _{1}\colon B] = 2$ and $-1 \not\in N(B _{1}/B)$.
Hence, $A _{-1} (-1, -1; B) \in d(B)$, and by Albert's theorem (cf.
[A2, Ch. IX, Sect. 6]), $B _{1}$ is not included in any cyclic
quartic extension of $E$. Note also that $A _{-1} (-1, -1; E ^{\prime
}) \not\in d(E ^{\prime })$, for any quadratic extension $E ^{\prime
}/E$. Since $A _{-1} (-1, -1; E ^{\prime })$ is $E ^{\prime
}$-isomorphic to $A _{-1} (-1, -1; E) \otimes _{E} E ^{\prime }$,
this is implied by (1.2) (ii) and the embeddability of $E ^{\prime }$
in $A _{-1} (-1, -1; E)$ over $E$. Therefore, $E ^{\prime }$ is a
nonreal field, so it follows from the Artin-Schreier theory (see [L1,
Ch. XI, Proposition 2]) that $E ^{\ast } = E ^{\ast 2} \cup -E ^{\ast
2}$ (where $-E ^{\ast 2} = \{-\lambda ^{2}\colon \ \lambda \in E
^{\ast }\}$). As $E$ is formally real, this proves that $E$ is
Pythagorean and $E _{1}$ is its unique quadratic extension in $E (2)$.
Hence, by Galois theory, finite proper extensions of $E$ in $E (2)$
are cyclic and include $E _{1}$ (see Remark 3.4 (ii)). Summing up the
obtained results, one concludes that $E (2) = E _{1}$, which
completes our proof.
\par
\medskip
The application of Lemma 3.5 offers no difficulty because a
Henselian valued field is formally real if and only if its residue
field is of the same kind (cf. [La, Theorem 3.16]).
\par
\medskip
{\bf Corollary 3.6.} {\it Let $(K, v)$ be a Henselian valued
stable field with $v(K)$ $p$-indivisible, for some prime $p$. Then
Br$(\widehat K) _{p}$ is divisible unless $p = 2$ and $\widehat K$
is formally real.}
\par
\medskip
{\it Proof.} If $p = {\rm char}(\widehat K)$ or $p \in {\cal
P}(\widehat K)$, this can be deduced from Witt's theorem (cf.
[Dr1, Sect. 15]), and from Proposition 2.1, Theorem 3.1 and Lemma
3.5, respectively. Also, it is known that if $\widehat K (2) = \widehat
K$, then Br$(\widehat K) _{2} = \{0\}$ (see [MS, (16.1)] and [A1,
Ch. VII, Sect. 9]). Suppose further that $p \not\in {\cal P}(\widehat
K)$, $p \ge 3$ and $p \neq {\rm char}(\widehat K)$, denote by
$\widehat K _{0}$ the prime subfield of $\widehat K$, and let
$\widehat K _{1}$ be the extension of $\widehat K _{0}$ generated
by the roots of unity in $\widehat K _{\rm sep}$ of $p$-primary
degrees. It is well-known that $\widehat K _{1} = \widehat K _{p}
(\varepsilon )$, where $\varepsilon $ is a primitive $p$-th root
of unity and $\widehat K _{p}$ is the unique $\hbox{\Bbb Z}
_{p}$-extension of $\widehat K _{0}$ in $\widehat K _{\rm sep}$.
At the same time, it is clear from Galois theory (cf. [Ko,
Proposition 2.11]) and the
assumptions on $p$ that $\widehat K _{p} \subseteq \widehat K$,
which implies that $\widehat K(\varepsilon ) = \widehat K\widehat
K _{1}$, i.e. $\widehat K(\varepsilon )$ contains primitive $p
^{n}$-th roots of unity, for all $n \in \hbox{\Bbb N}$. This
ensures that Br$(\widehat K(\varepsilon )) _{p}$ is divisible (see
[MS, (16.1)] and [P, Sect. 15.1, Corollary b]). Since $[\widehat
K(\varepsilon )\colon \widehat K]$ divides $p -
1$, one also concludes that the composition cor$_{\widehat
K(\varepsilon )/\widehat K} \circ \pi _{\widehat K/\widehat
K(\varepsilon )}$ induces an automorphism of Br$(\widehat K) _{p}$
(cf. [T, Theorem 2.5]). Now the divisibility of Br$(\widehat K)
_{p}$ is obtained in the same way as the proof of [FSS,
Proposition 2], so Corollary 3.6 is proved.
\par
\medskip
{\bf Remark 3.7.} Let $T$ be an abelian torsion group with divisible
$p$-components, for all prime numbers $p > 2$. It has been proved in
[Ch7] and [Ch8] that if the $2$-component of $T$ is divisible or of
order $2$, then $T$ is isomorphic to the Brauer group of the residue
field of some stable field $F(T)$ with a Henselian discrete
valuation (and in the former case, $F(T)$ can be chosen from the
class of absolutely stable fields). As the structure of divisible
abelian groups is known (cf. [F, Theorem 23.1]), Corollary 3.6,
Lemma 3.5 and this result fully describe the abelian torsion groups
realizable as Brauer groups of residue fields of Henselian valued
stable fields with totally indivisible value groups.
\par
\medskip
Lemma 3.5 and our next lemma provide a Galois-theoretic
characterization of the $p$-quasilocal property in the class of
fields with primitive $p$-th roots of unity:
\par
\medskip
{\bf Lemma 3.8.} {\it Let $E$ be a field containing a primitive
$p$-th root of unity, for some $p \in {\cal P}(E)$. Then $E$ is nonreal
and $p$-quasilocal if and only if ${\cal G}(E (p)/E)$ is a $p$-group of
Demushkin type.}
\par
\medskip
{\it Proof.} As noted, for example, in [Wa2], it follows from Galois
cohomology that there is a group isomorphism $\kappa \colon \ _{p}
{\rm Br}(E) \to H ^{2}({\cal G}(E(p)/E), \hbox{\Bbb F} _{p})$, such
that the diagram
\vskip0.5truecm \noindent (3.2) $\matrix{E ^{\ast }/E ^{\ast p}
\times E ^{\ast }/E ^{\ast p}
&\mapright{{\rm Symb}} &{}_p {\rm Br}(E) \cr
\mapdown{\mu \times \mu} && \mapdown{\kappa} \cr H ^{1} ({\cal G}(E
(p)/E), \hbox{\Bbb F} _{p}) \times H ^{1} ({\cal G}(E (p)/E),
\hbox{\Bbb F} _{p})
&\mapright{\cup } &H ^{2} ({\cal G}(E (p)/E), \hbox{\Bbb F} _{p})
\cr }$

\vskip0.5truecm \noindent is commutative, where $E ^{\ast p} = \{e
^{p}\colon \ e \in E ^{\ast }\}$, $\mu $ is the Kummer
isomorphism of $E ^{\ast }/E ^{\ast p}$ on $H ^{1} ({\cal G}(E (p)/E),
\hbox{\Bbb F} _{p})$, $\cup $ is the cup-product mapping and Symb
maps $(\alpha E ^{\ast p}, \beta E ^{\ast p})$ into $[A
_{\varepsilon } (\alpha , \beta ; E)]$, for every pair of
elements $\alpha , \beta \in E ^{\ast }$. This, combined with [MS,
(16.1)] and Lemmas 3.2 and 3.5, proves Lemma 3.8.
\par
\medskip
{\bf Remark 3.9.} Let $P$ be a nontrivial pro-$p$-group, for some
prime number $p$:
\par
(i) It is known that $P$ is a $p$-group of Demushkin type of degree
zero if and only if it is a free pro-$p$-group. When this is the
case, $P$ is determined, up-to an isomorphism, by its rank (cf.
[Se1, Ch. I, 4.2]), and is realizable as an absolute Galois group
of a field of any prescribed characteristic [LvdD, (4.8)] (see
also [Ch2, Remark 2.6]);
\par
(ii) We say that $P$ is a Demushkin group, if it is a $p$-group of
Demushkin type of degree one. The classification, up-to
isomorphisms, of the pro-$p$-groups of this kind and of finite or
countable rank is known (see [D1,2; Lab1,2] and [Se2]). We refer
the reader to [MW1,2], for a similar description of Demushkin
groups of countable ranks, which are realizable as Galois groups
of maximal $p$-extensions.
\par
\vskip0.75truecm
\centerline {\bf 4. The main result}
\par
\medskip
The main purpose of this paper is to prove the following theorem:
\par
\medskip
{\bf Theorem 4.1.} {\it Assume that $E$ is a $p$-quasilocal field
with respect to a prime number $p$, $R$ is an extension of $E$ in
$E (p)$, and $D \in d(E)$ is an algebra of $p$-primary index.
Then:}
\par
(i) {\it The homomorphism $\pi _{E/R}$ maps ${\rm Br}(E) _{p}$
surjectively on ${\rm Br}(R) _{p}$;}
\par
(ii) $R$ {\it is a $p$-quasilocal field;}
\par
(iii) {\it $R$ embeds in $D$ as an $E$-subalgebra if and only if
$[R\colon E]$ divides ${\rm ind}(D)$; $R$ is a splitting field of
$D$ if and only if $[R\colon E]$ is infinite or divisible by ${\rm
ind}(D)$;}
\par
(iv) {\it If $[R\colon E]$ is infinite, then ${\rm Br}(R) _{p} =
\{0\}$.}
\par
\medskip
Theorem 4.1 is proved in Section 7 on the basis of the following
two lemmas.
\par
\medskip
{\bf Lemma 4.2.} {\it Let $U$ be a cyclic extension of a field $E$
in $E (p)$, such that $[U\colon E] = p ^{t}$, and suppose that
$U$ splits each $\Omega \in d(E)$ whose exponent divides $p ^{t}$.
Then:}
\par
(i) $\pi _{E/U}$ {\it maps ${\rm Br}(E) _{p}$ surjectively on {\rm
Br}$(U) _{p}$ and ${\rm cor}_{U/E}$ maps {\rm Br}$(U) _{p}$
injectively into {\rm Br}$(E) _{p}$;}
\par
(ii) {\it If $\ {\rm Br}(E) _{p} = \{0\}$, then ${\rm Br}(U) _{p} =
\{0\}$ and $N(\Phi /U) = U ^{\ast }$, for every finite extension
$\Phi $ of $U$ in $E (p)$;}
\par
(iii) {\it If $\ \mu \in U ^{\ast }$ and $U _{1}$ is a cyclic
extension of $E$ in $E (p)$, such that $U _{1} \cap U = E$, then
$\mu \in N((U _{1}U)/U)$ if and only if $N _{E} ^{U} (\mu ) \in N(U
_{1}/E)$.}
\par
\medskip
{\it Proof}. (i) Recall first that Br$(U) _{p}$ is a module over the
integral group ring $\hbox{\Bbb Z} [{\cal G}(U/E)]$ with respect to
the group operation in Br$(U) _{p}$ and the multiplication from
$\hbox{\Bbb Z} [{\cal G}(U/E)] \times {\rm Br}(U) _{p}$ into Br$(U) _{p}$
induced canonically by the action of ${\cal G}(U/E)$ on $U$. Let $\psi $ be
a generator of ${\cal G}(U/E)$. By (1.5), the image of Br$(E) _{p}$ under
$\pi _{E/U}$ is the set $\{\theta \in {\rm Br}(U) _{p}\colon \ (\psi
- 1)\theta = 0\}$. Thus, it suffices for our first assertion to prove
that $\psi $ acts
trivially on Br$(U) _{p}$. For this, take any nonzero $\Delta \in
{\rm Br}(U) _{p}$; say exp$(\Delta) = p ^{s}$. The set $_{p^{s}}
{\rm Br}(U) = \{\eta \in {\rm Br}(U)\colon \ p ^{s}\eta = 0\}$ is a
$\hbox{\Bbb Z} [{\cal G}(U/E)]$-submodule of Br$(U) _{p}$ and can be
viewed as a module over the group ring $(\hbox{\Bbb Z}/p^{s}
\hbox{\Bbb Z}) [{\cal G}(U/E)]$. By Lemma 1.2, $\psi - 1$ is nilpotent in
$(\hbox{\Bbb Z}/p ^{s}\hbox{\Bbb Z}) [{\cal G}(U/E)]$, i.e. there is $k \in
\hbox{\Bbb N}$ with $(\psi - 1) ^{k}\Delta = 0$ in Br$(U) _{p}$.
Take $k$ minimal with this property. If $k = 1$, then $\psi \Delta =
\Delta$, as desired. So assume $k \ge 2$. Let $\Gamma = (\psi - 1)
^{k-2} \Delta $ (so $\Gamma = \Delta $ if $k=2$), and let $\Lambda
= (\psi - 1)\Gamma = (\psi - 1) ^{k-1}\Delta \neq 0$. Because
$(\psi - 1)\Lambda = (\psi - 1) ^{k}\Delta = 0$, there is $D \in
{\rm Br}(E) _{p}$ with $\pi _{E/U} (D) = \Lambda $ in Br$(U) _{p}$.
By [T, Theorem 2.5], then, in Br$(E)$,
$$p ^{t}D = {\rm cor}_{U/E} (\pi _{E/U} (D)) = {\rm cor}_{U/E}
(\Lambda ) =$$
$${\rm cor}_{U/E} (\psi \Gamma - \Gamma ) = {\rm cor}_{U/E}
(\psi \Gamma ) - {\rm cor}_{U/E} (\Gamma ) = 0.$$ Thus the
assumption on $U$ implies that $\pi _{E/U} (D) = 0$, i.e. $\Lambda =
(\psi - 1) ^{k-1}\Delta = 0$ in Br$(U) _{p}$, contradicting the
minimality of $k$. Hence, $k = 1$, and the first part of (i) is
proved.
\par
For the second assertion, take any $A$ in Br$(U) _{p}$ with
cor$_{U/E} (A) = 0$ (in Br$(E)$). We have just proved that $A = \pi
_{E/U} (B)$, for some $B \in {\rm Br}(E) _{p}$. Hence, by [T,
Theorem 2.5], $p ^{t}B = 0$, and by hypothesis, $A = \pi _{E/U} (B)
= 0$, proving the desired injectivity.
\par
(ii) The equality Br$(U) _{p} = \{0\}$ follows from the
inclusion Br$(U) _{p} \subseteq {\rm Im}(\pi _{E/U})$ and the
assumption that Br$(E) _{p} = \{0\}$. Note also that by Galois
theory and the subnormality of proper subgroups of finite $p$-groups
(see [L1, Ch. I, Sect. 6; Ch. VIII]), if $\Phi \neq U$, then $U$
has a proper cyclic extension $\Phi _{0}$ in $\Phi $. Therefore,
the statement that $N(\Phi /U) = U ^{\ast }$ can be deduced from
the triviality of Br$(U) _{p}$ and Br$(\Phi _{0}) _{p}$ by a
standard inductive argument relying upon (1.4) (ii) and the
transitivity of norm mappings in towers of finite extensions.
\par
(iii) The condition $U _{1} \cap U = E$ shows that if $\varphi $ is
an $E$-automorphism of $U _{1}$ of order $[U _{1}\colon E]$, then it
is uniquely extendable to a $U$-automorphism $\tilde \varphi $ of $U
_{1}U$ of the same order. Observing also that cor$_{U/E}$ maps $[((U
_{1}U)/U, \tilde \varphi , \mu )]$ into $[(U _{1}/E, \varphi , N
_{E} ^{U} (\mu ))]$ (e.g. as in the proof of [Ch6, (4.1) (iii)]),
one reduces Lemma 4.2 (iii) to a consequence of Lemma 4.2 (i).
\par
\medskip
We recall that the following lemma is proved in Sections 5 and 6.
\par
\medskip
{\bf Lemma 4.3.} {\it Let $E$ be a $p$-quasilocal field with
respect to a prime number $p$, and let $F _{1}$ and $F _{2}$ be
different extensions of $E$ in $E (p)$ of degree $p$. Then $N(F
_{1}/E)N(F _{2}/E) = E ^{\ast }$.}
\par
\medskip
Lemmas 4.2 and 4.3 enable one not only to establish the main
result of this paper but also to take a serious step towards
determining and characterizing the basic types of fields whose
finite abelian extensions and norm groups are related essentially
in the same way as in the classical local class field theory (see
[Ch6]). As to our next result, it supplements Theorem 4.1 by
showing that the class of $p$-quasilocal fields is closed under the
formation of purely inseparable extensions.
\par
\medskip
{\bf Proposition 4.4.} {\it Let $E$ be a field, $K/E$ a finite
purely inseparable extension and $p$ a prime number. Then $K$ is
$p$-quasilocal if and only if $E$ is $p$-quasilocal.}
\par
\medskip
{\it Proof.} Let $K _{1}$ be an extension of $K$ in $K _{\rm sep}$
and $E _{1}$ the separable closure of $E$ in $K _{1}$. Then $E
_{1}K = K _{1}$ and $E _{1} \cap K = E$ (see [L1, Ch. VII,
Sects. 4 and 7]). Hence, by Galois theory and [P, Sect. 9.2,
Proposition c], $K _{1}$ and $E _{1} \otimes _{E} K$ are
$K$-isomorphic. Note further that $K _{1}/K$ is Galois if and only
if $E _{1}/E$ is of the same kind; such being the case, the Galois
groups ${\cal G}(E _{1}/E)$ and ${\cal G}(K _{1}/K)$ are isomorphic (cf. [L1, Ch.
VIII, Theorem 4] and [Ko, Ch. 2]). It follows from (1.1), (1.2) and
the equality $K _{1} = E _{1}K$ that if $p \neq {\rm char}(E)$, $K
_{1}/K$ is a finite $p$-extension, $T \in d(E)$ and $[T] \in {\rm
Br}(E) _{p}$, then $K _{1}$ embeds in $T \otimes _{E} K$ as a
$K$-subalgebra if and only if $E _{1}$ embeds in $T$ over $E$.
Observing also that $\pi _{E/K}$ is surjective (by the
Albert-Hochschild theorem, see [Dr1, page 110] or [Se1, Ch. II,
2.2]), one concludes that it induces an isomorphism Br$(E) _{p}
\cong {\rm Br}(K) _{p}$. This completes the proof of Proposition
4.4 in the case of $p \neq {\rm char}(E)$. Henceforth, we assume
that $p = {\rm char}(E)$ and $K _{1}/K$ is cyclic of degree $p$.
Then $[K\colon E]$ is a power of $p$ and one may consider only the
special case where $[K\colon E] = p$ (cf. [L1, Ch. VII, Sect. 7]).
Our argument also relies on the existence (see [Dr1, Sect. 15,
Lemma 2]) of a cyclic extension $K _{2}$ of $K$ in $K (p)$, such
that $[K _{2}\colon K] = p ^{2}$ and $K _{1} \subset K _{2}$. Fix
a generator $\tau _{2}$ of ${\cal G}(K _{2}/K)$, put $E _{2} = K _{2}
\cap E (p)$ and let $\tau _{1}$, $\sigma _{1}$ and $\sigma _{2}$
be the automorphisms induced by $\tau _{2}$ on $K _{1}$, $E _{1}$
and $E _{2}$, respectively. We show that Proposition 4.4 can be
deduced from the following statements:
\par
\medskip
(4.1) (i) If $E$ is $p$-quasilocal, $D \in d(K)$ and ind$(D) = p$,
then $D$ is similar to the $K$-algebra $(K _{2}/K, \tau _{2}, d)$,
for some $d \in K ^{\ast }$;
\par
(ii) If $K$ is $p$-quasilocal, $\Delta \in d(E)$ and ind$(\Delta ) =
p$, then $\Delta $ is similar over $E$ to $(E _{2}/E, \sigma _{2},
\delta )$, for some $\delta \in E ^{\ast }$;
\par
(iii) If $E$ or $K$ is $p$-quasilocal, then $N(E _{2}/E _{1})$ and
$N(K _{2}/K _{1})$ include the sets
\par \noindent
$\{\sigma _{1} (\alpha )\alpha ^{-1}\colon \ \alpha \in E _{1}
^{\ast }\}$ and $\{\tau _{1} (\beta )\beta ^{-1}\colon \ \beta \in
K _{1} ^{\ast }\}$, respectively.
\par
\medskip
Suppose first that $E$ is $p$-quasilocal. Then it follows from (1.4)
(ii) and (4.1) (i) that $d ^{p} \in N(K _{2}/K)$, and by Hilbert's
Theorem 90, this means that $d\tau _{1} (\alpha _{1})\alpha _{1}
^{-1} \in N(K _{2}/K _{1})$, for some $\alpha _{1} \in K _{1} ^{\ast
}$. Applying now (4.1) (iii), one concludes that $d \in N(K _{2}/K
_{1})$. The obtained result is equivalent to the embeddability of $K
_{1}$ in $D$ as a $K$-subalgebra, since $(K _{2}/K, \tau _{2}, d)
\otimes _{K} K _{1}$ is similar over $K _{1}$ to the centralizer of
$K _{1}$ in $(K _{2}/K, \tau _{2}, d)$, i.e. to $(K _{2}/K _{1},
\tau _{2} ^{p}, d)$ (cf. [P, Sect. 13.3]). This proves that $K$ is
$p$-quasilocal. The converse implication of Proposition 4.4 is
deduced from (4.1) (ii)-(iii) in much the same way, so we omit the
details.
\par
We turn to the proof of (4.1). Denote by $\pi _{F/\widetilde
F,p}$ the scalar extension map of Br$(F) _{p}$ into Br$(\widetilde
F) _{p}$, for each pair $(F, \widetilde F)$ of intermediate fields
of $K (p)/E$ satisfying the inclusion $F \subseteq \widetilde F$.
It is not difficult to see from (1.4), (1.5) and Lemma 4.2 (i)
that (4.1) (iii) will be proved, if we show that $\pi _{E/E
_{1},p}$ is surjective if and only if so is $\pi _{K/K _{1},p}$. Our
proof of this equivalence relies on the fact that $\pi _{E/K
_{1},p}$ equals the compositions $\pi _{E _{1}/K _{1},p} \circ \pi
_{E/E _{1},p}$ and $\pi _{K/K _{1},p} \circ \pi _{E/K,p}$ (see [P,
Sect. 9.4, Corollary a]). Since $K/E$ and $K _{1}/E _{1}$ are purely
inseparable, $\pi _{E/K,p}$ and $\pi _{E _{1}/K _{1},p}$ are
surjective, so one obtains consecutively that if $\pi _{E/E _{1},p}$
is surjective, then the same applies to $\pi _{E/K _{1},p}$ and $\pi
_{K/K _{1},p}$. Conversely, let $\pi _{K/K _{1},p}$ be surjective.
Then $\pi _{E/K _{1},p}$ is surjective, which implies that Br$(E
_{1}) _{p}$ is presentable as a sum of the subgroups Im$(\pi _{E/E
_{1},p})$ and Br$(K _{1}/E _{1})$. Since Br$(E _{1}) _{p}$ and Im$(\pi
_{E/E _{1},p})$ are divisible and Br$(K _{1}/E _{1})$ is of exponent
dividing $[K _{1}\colon E _{1}] = [K\colon E]$, the obtained result
proves the surjectivity of $\pi _{E/E _{1},p}$ (and the validity of
(4.1) (iii)).
\par
For the rest of the proof, note that if $E$ is $p$-quasilocal and
$D _{1} \in d(E)$ is chosen so that $\pi _{E/K} ([D _{1}]) = [D]$,
then exp$(D _{1})$ divides $p.[K\colon E] = p ^{2} = [E _{2}\colon
E]$ (apply (1.1) (i) and (1.2) (ii)). Hence, by Theorem 4.1 (iii),
$D _{1}$ is split by $E _{2}$, which allows one to deduce (4.1)
(i) from (1.4) and the existence of a $K$-isomorphism $K _{2}
\cong E _{2} \otimes _{E} K$. Let now $K$ be $p$-quasilocal. Then
$K _{2}$ splits all algebras in $s(K)$ of exponents dividing $p
^{2}$. In particular, this applies to the $K$-algebra $\Theta
\otimes _{E} K$ whenever $\Theta \in d(E)$ and exp$(\Theta ) = p
^{2}$. Since Br$(E) _{p}$ is divisible, $\pi _{E/K _{2},p} = \pi
_{E _{2}/K _{2},p} \circ \pi _{E/E _{2},p}$ and $[K _{2}\colon E
_{2}] = [K\colon E] = p$, this implies that $_{p} {\rm Br}(E)
\subseteq {\rm Br}(E _{2}/E)$ and so completes the proof of (4.1)
(ii) and Proposition 4.4.
\par
\medskip
Note that Theorem 3.1 enables one to verify more easily whether a
given Henselian valued field is stable. This can be illustrated by
the following result:
\par
\medskip
{\bf Proposition 4.5.} {\it Assume that $(K, v)$ is a Henselian
valued field such that $v(K)$ is totally indivisible, ${\rm
char}(\widehat K) = q \ge 0$ and ${\cal P}(\widehat K)$ contains every
prime $p$ for which {\rm Br}$(\widehat K) _{p} \neq \{0\}$.
Suppose also that if {\rm Br}$(K) _{p'} \neq \{0\}$, for some
prime $p ^{\prime }$, then $_{p'} {\rm Br}(K)$ coincides with the
set $\{[D _{p'}]\colon \ D _{p'} \in d(K), {\rm ind}(D _{p'}) = p
^{\prime }\}$. Then:}
\par
(i) ind$(D) = {\rm exp}(D)${\it , for every $D \in d(K)$ with
$[D\colon K]$ not divisible by ${\rm char}(\widehat K)$;}
\par
(ii) {\it $K$ is stable in each of the following three
special cases:}
\par
($\alpha $) $K$ {\it is almost perfect and ${\rm char}(K) = q$;}
\par
($\beta $) $\widehat K$ {\it is perfect, $q > 0$,
the group $v(K)/qv(K)$ is of order $q$ and $d(K)$ is included in
the class of defectless division $K$-algebras;}
\par
($\gamma $) $K$ {\it contains a primitive $q$-th root of unity and
$v(K)/qv(K)$ is of order $q ^{3}$.}
\par
\medskip
It is not known whether there exists a field $E$ and a prime $p
\not\in {\cal P}(E)$, for which Br$(E) _{p} \neq \{0\}$. In view of the
Merkurjev-Suslin theorem and [A1, Ch. VII, Theorem 28], this
is impossible, if $E$ contains a primitive $p$-th root of unity or
char$(E) = p$.
\par
\medskip
{\it Proof of Proposition 4.5.} Let $\widetilde \Delta \in
d(\widehat K)$ be of prime index $p$, and suppose that $\widetilde
L/\widehat K$ is a cyclic extension of degree $p$, $\Delta $ and $L$
are inertial lifts over $K$ of $\widetilde \Delta $ and $\widetilde
L$, respectively, $\sigma $ is a generator of ${\cal G}(L/K)$ (see [JW,
page 135]), and $\pi $ is an element of $K ^{\ast }$ of value $v(\pi
) \not \in pv(K)$. As in the proof of Proposition 2.1, one deduces
from [JW, Theorem 5.15 (a)] that if $\widetilde L$ does not embed in
$\widetilde \Delta $ as a $\widehat K$-subalgebra, then $\Delta
\otimes _{K} (L/K, \sigma , \pi )$ must lie in $d(K)$ and have
exponent $p$ and index $p ^{2}$. This contradicts the assumptions of
Proposition 4.5 and so proves that $\widehat K$ is $p$-quasilocal.
Hence, by Theorems 3.1, 4.1 and the condition on ${\cal P}(\widehat
K)$, $\widehat K$ is stable and its cyclic extensions have the property
required by Corollary 2.2. Observe also that our assumptions rule
out the existence in $d(K)$ of a tensor product of a pair of cyclic
division $K$-algebras of any prime index $p ^{\prime }$. As shown in
[Ch1, Sect. 2], this means that $K, \widehat K, \widehat K (p
^{\prime }), v(K)$ and $v(K)/p ^{\prime }v(K)$ are related as in
[Ch1, Theorem 2.1]. Therefore, our proof can be completed by
applying [Ch1, Theorem 3.1].
\par
\medskip
We conclude this Section with essentially an equivalent form of
Theorem 4.1 for nonreal fields containing a primitive $p$-th root
of unity. It partially generalizes Theorem 2 of [Lab1]:
\par
\medskip
{\bf Proposition 4.6.} {\it Let $E$ be a field containing a
primitive $p$-th root of unity, for some prime $p$, and let ${\cal
G}(E (p)/E)$ be a $p$-group of Demushkin type of degree $d$. Then
${\cal G}(E (p)/E)$ has the following properties:}
\par
(i) {\it Every open subgroup $U$ of ${\cal G}(E (p)/E)$ is a $p$-group of
Demushkin type of degree $d$; the corestriction mapping of $H ^{2}
(U, \hbox{\Bbb F} _{p})$ into $H ^{2} ({\cal G}(E (p)/E), \hbox{\Bbb F}
_{p})$ is an isomorphism;}
\par
(ii) {\it Nontrivial closed subgroups of ${\cal G}(E (p)/E)$ of infinite
indices are free pro-$p$-groups.}
\par
\medskip
{\it Proof.} Since $E (p) = E ^{\prime } (p)$, for every extension
$E ^{\prime }$ of $E$ in $E (p)$ (cf. [Ko, the beginning of Ch. 9]),
this can be deduced from Theorem 4.1, Lemmas 3.8 and 4.2 (i),
diagram (3.2) and [Se1, Ch. I, Proposition 21].
\par
\medskip
{\bf Remark 4.7.} It is likely that Proposition 4.6 (i) remains
valid for a large class of pro-$p$-groups, including Poincar$\acute
e$ groups of arbitrary dimensions and $p$-groups of Demushkin type
(see, for example, [Se1, Ch. I, 4.5]).
\par
\medskip
It is easily seen that the degree of any $p$-group of Demushkin type
is at most equal to its rank as a pro-$p$-group. Note also that if
$F$ is a field with a primitive $p$-th root of unity and ${\cal G}(F
(p)/F)$ of Demushkin type of degree $\ge 2$, then the rank of ${\cal
G}(F (p)/F)$ is infinite [Ch6, Corollary 4.6]. It is therefore worth
mentioning (for a proof, see [Ch8]) that if $E _{0}$ is an infinite
field of cardinality $d$ and characteristic $q \ge 0$, ${\cal P}
_{q}$ is the set of prime numbers different from $q$, and $c _{p}, d
_{p}\colon \ p \in {\cal P} _{q}$, is a system of cardinal numbers
such that $c _{p} \le d _{p} \le d, \ p \in {\cal P} _{q}$, then
there exists an extension $E$ of $E _{0}$ with the following
properties:
\par
\medskip
(4.2) $E _{0}$ is algebraically closed in $E$, and for each $p \in
{\cal P} _{q}$, ${\cal G}(E (p)/E)$ and the Sylow pro-$p$-subgroups
of ${\cal G} _{E}$ are of rank $d$ and Demushkin type of degrees $c
_{p}$ and $d _{p}$, respectively. Moreover, $E$ can be chosen so
that finite groups are realizable as Galois groups over $E$ (and $E$
does not admit Henselian valuations with indivisible value groups,
cf. [Ch6, (2.3)]).
\par
\medskip
This, applied to the special case where $q = 0$, $d = \aleph _{0}$
and $c _{p} = d _{p} = 1$, for every $p \in {\cal P} _{q}$, enables one
to deduce from the main results of [MW1,2] that each sequence $G
_{p}\colon \ p \in {\cal P} _{q}$, of Demushkin groups of countable rank
and $s$-invariant zero is realizable as a sequence of Sylow
pro-$p$-subgroups of the absolute Galois group of some field.
\par
\vskip0.75truecm
{\bf 5. An equivalent form of the main lemma}
\par
\medskip
Our aim in this Section is to find an equivalent form of Lemma 4.3.
This form is of independent interest and serves as a basis for the
proof of the lemma, presented in Section 6. The main result of this
Section is known in the special case of $p = 2$ and a ground field
of characteristic different from $2$ (see Exercise 4.4 at the end of
[CF]).
\par
\medskip
{\bf Lemma 5.1.} {\it Let $E$ be a field, $L$ an extension of $E$
of degree $p ^{2}$, for some prime $p$, and $I(L/E)$ the set of
intermediate fields of $L/E$. Then $L/E$ is noncyclic and Galois
if and only if it satisfies any of the following equivalent
conditions:}
\par
(i) {\it $L$ is a compositum of two different cyclic extensions
of $E$ of degree $p$;}
\par
(ii) $I(L/E) = \{L, E, E _{j}\colon \ j = 1, \dots , p + 1\}${\it,
where $E _{1}, \dots , E _{p+1}$ are (pairwise distinct) cyclic
extensions of $E$ of degree $p$.}
\par
\medskip
{\it Proof.} This follows at once from Galois theory and the
well-known fact that noncyclic groups of order $p ^{2}$ are
elementary abelian with exactly $p + 1$ subgroups of order $p$.
\par
\medskip
{\bf Lemma 5.2.} {\it Let $E$ be a field containing a primitive
$p$-th root of unity $\varepsilon $, $L/E$ a noncyclic Galois
extension of degree $p ^{2}$, and $(F _{1}, F _{2})$ a pair of
different extensions of $E$ in $L$ of degree $p$. Then $F _{i} =
E(\xi _{i})$, where $\xi _{i}$ is a $p$-th root of an element $a
_{i} \in E ^{\ast }$, for each index $i$. Moreover,}
\par
(i) {\it The multiplicative group $F _{i} ^{\ast p} \cap E$ equals
the union $\bigcup _{l=0} ^{p-1} a _{i} ^{l}.E ^{\ast p}$, and $L
^{\ast p} \cap E =$ $(F _{1} ^{\ast p} \cap E).(F _{2} ^{\ast p}
\cap E)$;}
\par
(ii) {\it If $p > 2$, then $N _{E} ^{F _{i}} (\xi _{i}) = a _{i}$
and $F _{i} ^{\ast p} \cap E$ is a subgroup of $N(F _{i}/E)$, for
$i = 1, 2$. Furthermore, one of the following conditions holds:}
\par
(a) $(L ^{\ast p} \cap N(F _{1}/E)) = (F _{1} ^{\ast p} \cap E)${\it
; this occurs if and only if $(L ^{\ast p} \cap N(F _{2}/E)) = (F
_{2} ^{\ast p}$ $\cap \ E)$;}
\par
(b) $(L ^{\ast p} \cap E) \subseteq N(F _{1}/E)${\it ; this is
the case if and only if $(L ^{\ast p} \cap E) \subseteq N(F
_{2}/E)$.}
\par
\medskip
{\it Proof.} The existence of $a _{i}, \xi _{i}\colon \ i = 1, 2$,
and statement (i) follow at once from Kummer theory, and the
former part of (ii) is implied by the definition of the norm
mapping. In view of (1.4) (ii), the obtained results indicate that
$(L ^{\ast p} \cap E) \subseteq N(F _{1}/E)$ if and only if $A
_{\varepsilon } (a _{1}, a _{2}; E) \not\in d(E)$. Similarly, we
have $(L ^{\ast p} \cap E) \subseteq N(F _{2}/E)$ if and only if
$A _{\varepsilon } (a _{2}, a _{1}; E) \not\in d(E)$. Since $A
_{\varepsilon } (a _{1}, a _{2}; E)$ and $A _{\varepsilon } (a
_{2}, a _{1}; E)$ are inversely isomorphic over $E$, this proves
Lemma 5.2.
\par
\medskip
{\bf Lemma 5.3.} {\it Let $E$ be a field not containing a
primitive $p$-th root of unity, for some prime $p$. Then $E (p)$
does not contain such a root and $E (p) ^{\ast p} \cap E = E
^{\ast p}$.}
\par
\medskip
{\it Proof.} If char$(E) = p$, this follows from the fact that the
binomial $X ^{p} - \alpha $ is purely inseparable, for each
$\alpha \in E$. Assuming that $p \neq {\rm char}(E)$ and $\alpha \in
E ^{\ast } \setminus E ^{\ast p}$, one obtains that $X ^{p} -
\alpha $ is irreducible over $E$ and its root field, say $F
_{\alpha }$, contains a primitive $p$-th root of unity $\varepsilon
$. As $[E(\varepsilon )\colon E]$ divides $p - 1$ and
$\varepsilon \not\in E$, this means that the extension of $E$
generated by a fixed $p$-th root of $\alpha $ in $F _{\alpha }$ is
not normal. Our conclusion, however, contradicts the normality of
the extensions of $E$ in $E (p)$ of degree $p$, so Lemma 5.3 is
proved.
\par
\medskip
{\bf Lemma 5.4.} {\it Let $E$ be a field, $p$ a prime number, $L/E$
a noncyclic Galois extension of degree $p ^{2}$, $E _{1}, \dots , E
_{p+1}$ the intermediate fields of $L/E$ of degree $p$ over $E$, and
$\sigma _{1}, \dots , \sigma _{p+1}$ generators of ${\cal G}(E _{1}/E),
\dots , {\cal G}(E _{p+1}/E)$, respectively. Then $N _{E} ^{L} (\alpha
)\alpha ^{p} = \prod _{j=1} ^{p+1} N _{E _{j}} ^{L} (\alpha )$, for
every $\alpha \in L$. Furthermore, if $N _{E} ^{L} (\alpha ) = c
^{p}$, for some $c \in E ^{\ast }$, then there exist elements $u
_{1} \in E _{1} ^{\ast }, \dots , u _{p+1} \in E _{p+1}
^{\ast }$, such that $N _{E _{j}} ^{L} (\alpha ) = c\sigma _{j} (u
_{j})u _{j} ^{-1}$, for every index $j$. In addition, $\alpha ^{p} =
c.\prod _{j=1} ^{p+1} (\sigma _{j} (u _{j})u _{j} ^{-1})$.}
\par
\medskip
{\it Proof.} This follows from Hilbert's Theorem 90 and the
definitions of the considered norm mappings.
\par
\medskip
From now on we will often use the fact that if $M/E$ is a Galois
extension, then $M ^{\ast }$ is a module over the integral group
ring $\hbox{\Bbb Z} [{\cal G}(M/E)]$ with respect to the group operation
in $M ^{\ast }$ and the multiplication $\ast \colon \hbox{\Bbb Z}
[{\cal G}(M/E)] \times M ^{\ast } \to M ^{\ast }$ canonically induced by
the action of ${\cal G}(M/E)$ on $M ^{\ast }$. Note also that if ${\cal
G}(M/E)$ is a finite abelian group and $F$ is an intermediate field of
$M/E$, then $F ^{\ast }$, $N(M/F)$ and $F ^{\ast l}$ are
$\hbox{\Bbb Z} [{\cal G}(M/E)]$-submodules of $M ^{\ast }$ satisfying the
inclusions $F ^{\ast l} \subseteq N(M/F) \subseteq F ^{\ast }$,
where $l = [F\colon E]$.
\par
\medskip
{\bf Lemma 5.5.} {\it Let $E$ be a field, $F$ a cyclic extension
of $E$ of prime degree $p$, $\alpha $ an element of $F ^{\ast }$,
$k$ an integer with $0 < k < p$, $\tau _{i}\colon \ i = 1, \dots ,
k$ a sequence of generators of ${\cal G}(F/E)$, and $T _{i} = \tau _{i} -
1 \in \hbox{\Bbb Z} [{\cal G}(F/E)]$, for each index $i$. Then:}
\par
(i) {\it Suppose $k = p - 1$, $\sigma $ is any generator of
${\cal G}(F/E)$ and each $\tau _{i} = \sigma ^{m(i)}$. Then $(\prod
_{i=1} ^{p-1} T _{i}) \ast \alpha = (N _{E} ^{F} (\alpha )\alpha
^{-p}) ^{m} (\sigma (\beta )\beta ^{-1}) ^{p}$ for some $\beta \in
F ^{\ast }$, where $m = \prod _{i=1} ^{p-1} m(i)$;}
\par
(ii) {\it If $(\prod _{i=1} ^{k} T _{i}) \ast \alpha = 1$, then
$\alpha ^{p} \in E$ (and $N _{E} ^{F} (\alpha ) = \alpha ^{p}$);
the converse is true in case $k \ge 2$;}
\par
(iii) {\it If $p \ge 3$, $k \ge 3$ and $(\prod _{i=1} ^{k-1} T _{i})
\ast \alpha = \rho $, with $N _{E} ^{F} (\alpha ) \in F ^{\ast p}$,
then
\par \noindent
$\rho = (\prod _{i=1} ^{k} T _{i}) \ast \gamma
$, for some $\gamma \in F ^{\ast }$.}
\par
\medskip
{\it Proof.} (i) In the polynomial ring $\hbox{\Bbb Z} [X]$ (in
one indeterminate), let $Y = X - 1$. Note first that for some
$g(X) \in \hbox{\Bbb Z} [X]$,
$$(5.1) \ - p + \sum _{i=0} ^{p-1} X ^{i} = (X - 1) ^{p-1} + p(X -
1)g(X).$$
This is easy to see from the binomial expansion in terms of $Y$:
$$ \ \sum _{i=0} ^{p-1} X ^{i} = (X ^{p} - 1)/(X - 1) = [(Y +
1) ^{p} - 1]/Y = Y ^{p-1} + p + \sum _{j=2} ^{p-1} \pmatrix{ p \cr j
\cr } Y ^{j-1}.$$ Now, note that for any $l \in \hbox{\Bbb N}$, we
have $X ^{l} - 1 = (X - 1)(X ^{l-1} + \dots + 1) = (X - 1)[(X$ $-
1)h _{l} (X) + l]$, for some $h _{l} (X) \in\hbox{\Bbb Z} [X]$.
Then, for some $h(X), f(X), q(X) \in \hbox{\Bbb Z} [X]$, using (5.1)
at the last but one step,
$$\prod _{i=1} ^{p-1} (X ^{m(i)} - 1) = (X - 1) ^{p-1} \prod
_{i=1} ^{p-1} [(X - 1)h _{m(i)} (X) + m(i)]$$
$$= (X - 1) ^{p-1} [m + (X - 1)h(X)]$$
$$= m(X - 1) ^{p-1} + [(X - 1) ^{p} - X ^{p} + 1)]h(X) +
(X ^{p} - 1)h(X)$$
$$= m[\sum _{i=0} ^{p-1} X ^{i} - p - p(X - 1)g(X)] + p(X - 1)f(X)
+ (X ^{p} - 1)h(X)$$
$$= m[\sum _{i=0} ^{p-1} X ^{i} - p] + p(X - 1)q(X) + (X ^{p} -
1)h(X).$$
The formula in Lemma 5.5 (i) is obtained by evaluating this
equation in the group ring $\hbox{\Bbb Z} [{\cal G}(F/E)]$, mapping $X$
into $\sigma $, then applying the result to $\alpha $.
Specifically, $\beta = q(\sigma ) \ast \alpha $.
\par
(ii) It is clear from Galois theory that an element $\alpha \in F
^{\ast }$ satisfies the equality $\tau _{1} (\alpha )\alpha ^{-1}
= 1$ if and only if $\alpha \in E$. When $\alpha \in E$, we have
$N _{E} ^{F} (\alpha ) = \alpha ^{p}$, so the former part of our
assertion is proved in the case of $k = 1$. The obtained result
also indicates that if $(T _{1}T _{2}) \ast \alpha = 1$, then
$\tau _{2} (\alpha )\alpha ^{-1}$ is a $p$-th root of unity lying
in $E$, and since $p > 2$, this yields $N _{E} ^{F} (\alpha ) =
\alpha ^{p}$. Suppose further that $k \ge 2$. Clearly, if $\alpha
^{p} \in E ^{\ast }$, then $\tau _{k} (\alpha )\alpha ^{-1}$ is a
$p$-th root of unity. In view of Lemma 5.3, this root lies in $E$,
so $(T _{k-1}T _{k}) \ast \alpha = 1$, which proves the latter
part of Lemma 5.5 (ii). It remains to be seen that $\alpha ^{p}
\in E$, provided that $k \ge 3$ and $(\prod _{i=1} ^{k} T _{i})
\ast \alpha = 1$. Then the element $(\prod _{i=1} ^{k-2} T _{i})
\ast \alpha := \bar \alpha $ satisfies the equality $N _{E} ^{F}
(\bar \alpha ) = \bar \alpha ^{p} = 1$. Moreover, by Lemma 5.3,
$\bar \alpha \in E ^{\ast }$, which implies that $(\prod _{i=1}
^{k-1} T _{i}) \ast \alpha = 1$. This result, used repeatedly,
leads to the conclusion that $(T _{1}T _{2}) \ast \alpha = 1$, and
so completes the proof of the former part of Lemma 5.5 (ii).
\par
(iii) Suppose that $N _{E} ^{F} (\alpha ) = \alpha _{0} ^{p}$, for
some $\alpha _{0} \in F ^{\ast }$. By Hilbert's Theorem 90, then
$\alpha = \alpha _{0}\tau _{k} (\gamma )\gamma ^{-1}$, for some
$\gamma \in F ^{\ast }$. Hence, the inequalities $p \ge 3$, $k \ge
3$ and the latter part of (ii), applied to $\alpha _{0}$, yield
$(\prod _{i=1} ^{k-1} T _{i}) \ast \alpha = (\prod _{i=1} ^{k} T
_{i}) \ast \gamma $, so Lemma 5.5 is proved.
\par
\medskip
The main result of this Section can be stated as follows:
\par
\medskip
{\bf Proposition 5.6.} {\it Assume that $p$ is a prime number, and
$L/E$ is a noncyclic Galois extension of fields with $[L\colon
E] = p ^{2}$. Let $E _{1}, \dots , E _{p+1}$ be the intermediate
fields, $E \subset E _{i} \subset L$, and for each $i \le p - 1$, let
$\varphi _{i}$ be any generator of ${\cal G}(L/E _{i})$ and $N _{i} =
\varphi _{i} - 1 \in \hbox{\Bbb Z} [{\cal G}(L/E _{i})]$. Then the
following conditions are equivalent, for any $c \in E ^{\ast }$:}
\par
(i) {\it $c \in N(E _{p}/E)N(E _{p+1}/E)$;}
\par
(ii) {\it There exist elements $\zeta \in L, z _{p} \in E
_{p}$ and $z _{p+1} \in E _{p+1}$ such that
$$N _{E _{i}} ^{L} (\zeta ) = c \ {\it for} \ 1 \le i \le p - 1, \
{\it and} \ N _{E _{j}} ^{L} (\zeta ) = c.(\prod _{t=1} ^{p-1} N
_{t}) \ast z _{j} \ {\it for} \ j = p, p + 1.$$}
\par
\medskip
{\it Proof.} For each pair of indices $j \ge p$, $i \le p - 1$,
$\varphi _{i}$ induces on $E _{j}$ an automorphism of order $p$,
so Lemma 5.5 (i) implies the following statement:
\par
\medskip
(5.2) There exists a positive integer $f(j)$ not divisible by $p$,
such that
\par \noindent
$(\prod _{t=1} ^{p-1} N _{t}) \ast (\mu
_{j})E _{j} ^{\ast p} = N _{E} ^{E _{j}} (\mu _{j} ^{f(j)})E _{j}
^{\ast p}$, for every $\mu _{j} \in E _{j} ^{\ast }$.
\par
\medskip
We first show that (ii)$\to $(i) Suppose that $\zeta , z _{p}$ and
$z _{p+1}$ are related as in Proposition 5.6 (ii). By Lemma 5.4,
then $\zeta ^{p} = c.(\prod _{t=1} ^{p-1} N _{t}) \ast (z _{p}z
_{p+1})$, so it follows from (5.2) that
\par \noindent
$cN _{E} ^{E _{p}} (z _{p} ^{-f(p)})N _{E} ^{E _{p+1}} (z _{p+1}
^{-f(p+1)})$ is contained in $L ^{\ast p}$. Applying now Lemmas
5.2 (ii) and 5.3, one concludes that $c \in N(E _{p}/E)N(E
_{p+1}/E)$ except, possibly, in the special case of $p = 2 \neq
{\rm char}(E)$ and $\sqrt{-1} \not\in E$. At the same time, it is
easily verified that if $p = 2$, then
$$\varphi _{1} (z _{3})z _{3} ^{-1} = (\varphi _{3}\varphi _{1})
(z _{3})z _{3} ^{-1} = N _{E _{3}} ^{L} (\zeta )N _{E _{1}} ^{L}
(\zeta ) ^{-1} = (\varphi _{3}\varphi _{1}) (\varphi _{1} (\zeta
))\varphi _{1} (\zeta ) ^{-1}$$ and $\varphi _{3}\varphi _{1} =
\varphi _{2}$. Hence, by Galois theory, $E _{2} ^{\ast }$ contains
the element $z _{2} ^{\prime } = z _{3} ^{-1}\varphi _{1} (\zeta )$.
Therefore, we have
$$c = N _{E _{1}} ^{L} (\zeta ) = N _{E _{1}} ^{L} (\varphi _{1}
(\zeta )) = N _{E _{1}} ^{L} (z _{3})N _{E _{1}} ^{L} (z _{2}
^{\prime }) = N _{E} ^{E _{3}} (z _{3})N _{E} ^{E _{2}} (z _{2}
^{\prime }),$$
which completes our proof.
\par
We prove that (i)$\to $(ii) Let $c = N _{E} ^{E _{p}} (\alpha _{p})N
_{E} ^{E _{p+1}} (\alpha _{p+1})$, for some $\alpha _{p} \in E
_{p}$, $\alpha _{p+1} \in E _{p+1}$. Since $N _{E _{j'}} ^{L}
(\alpha _{j}) = N _{E} ^{E _{j}} (\alpha _{j})$, for each $j \ge p$
and any index $j ^{\prime } \neq j$, the product $\zeta = \alpha
_{p}\alpha _{p+1}$ satisfies the equalities $N _{E _{i}} ^{L} (\zeta
) = c\colon \ i \le p - 1$, and $N _{E _{j}} ^{L} (\zeta ) = ca
_{j}\colon \ j = p$, $p + 1$, where $a _{j} = \alpha _{j} ^{p}N _{E}
^{E _{j}} (\alpha _{j}) ^{-1}$. Therefore, (i)$\to $(ii) will be
proved, if we show that the equation $(\prod _{t=1} ^{p-1} N _{t})
\ast X _{j} = a _{j}$ has a solution in $E _{j} ^{\ast }$, for each
$j \ge p$. Clearly, $N _{E} ^{E _{j}} (a _{j}) = 1$, and by
Hilbert's Theorem 90, this yields $a _{j} = \varphi _{1} (b _{j})b
_{j} ^{-1}$, for some $b _{j} \in E _{j} ^{\ast }$. When $p = 2$,
our assertion is thereby proved, so we assume further that $p > 2$.
Fix an integer $k _{j}$ so that $k _{j}m _{j} \equiv - 1 ({\rm mod}
p)$, where $m _{j}$ is determined as $m$ in Lemma 5.5 (i) by the
restrictions of $\varphi _{1}$ and $\varphi _{t}\colon \ t = 1,
\dots , p - 1$, on $E _{j}$. Applying Lemma 5.5 (i) and the equality
$a _{j} = \varphi _{1} (b _{j})b _{j} ^{-1}$, one obtains that $a
_{j} = [(\prod _{t=1} ^{p-1} N _{t}) \ast \alpha _{j} ^{k
_{j}}]\varphi _{1} (\gamma _{j} ^{p})\gamma _{j} ^{-p}$, for some
$\gamma _{j} \in E _{j} ^{\ast }$. Now it suffices for the proof of
(i)$\to $(ii) to establish the solvability of the equation $(\prod
_{t=1} ^{p} N _{t}) \ast X _{j} = \varphi _{1} (\gamma _{j})
^{p}\gamma _{j} ^{-p}$ over $E _{j} ^{\ast }$, where $N _{p} = N
_{1}$. This can be stated more completely as follows:
\par
\medskip
(5.3) For an element $\rho _{j}$ of $E _{j} ^{\ast }$, the
following conditions are equivalent:
\par
(c) The equation $(\prod _{t=1} ^{p} N _{t}) \ast X _{j} = \rho
_{j}$ is solvable over $E _{j}$;
\par
(cc) $\rho _{j} \in E _{j} ^{\ast p}$ and the equation $(N _{1}N
_{2}) \ast Y _{j} = \rho _{j}$ is solvable over $E _{j}$;
\par
(ccc) There exists an element $\eta _{j} \in E _{j} ^{\ast }$ for
which $\varphi _{1} (\eta _{j}) ^{p}\eta _{j} ^{-p} = \rho _{j}$.
\par
\medskip
We prove (5.3). (c)$\to $(ccc) Let $\tilde \rho _{j}$ be an
element of $E _{j}$, such that $(\prod _{t=1} ^{p} N _{t}) \ast
\tilde \rho _{j} = \rho _{j}$. By (5.2), we have $(\prod
_{t=1} ^{p-1} N _{t}) \ast \tilde \rho _{j} = N _{E} ^{E _{j}}
(\tilde \rho _{j} ^{f(j)})\eta _{j} ^{p}$, for some $f(j) \in
\hbox{\Bbb Z}$ and $\eta _{j} \in E _{j}$. It is therefore clear
that $\rho _{j} = \varphi _{1} (\eta _{j}) ^{p}\eta _{j} ^{-p}$,
whence (c)$\to $(ccc).
\par
(cc)$\to $(c) Assume that $\rho _{j} \in E _{j} ^{\ast p}$ and
$(N _{1}N _{2}) \ast \rho _{j} ^{\prime } = \rho _{j}$, for some
$\rho _{j} ^{\prime } \in E _{j}$. We prove (c) by assuming the
opposite. Then one obtains, using repeatedly Lemma 5.5 (iii),
that there exists a pair $(n(j), \bar \rho _{j}) \in (\hbox{\Bbb
Z} \times E _{j} ^{\ast })$, such that $2 \le n(j) \le p - 1$,
$(\prod _{t=1} ^{n(j)} N _{t}) \ast \bar \rho _{j} = \rho _{j}$
and $N _{E} ^{E _{j}} (\bar \rho _{j}) \not\in E _{j} ^{\ast p}$.
On the other hand, Lemma 5.5 (i) indicates that if $n(j) < p - 1$,
then $(\prod _{t=n(j)+1} ^{p-1} N _{t}) \ast \rho _{j}$ could not
lie in $E _{j} ^{\ast p}$, which contradicts the condition $\rho
_{j} \in E _{j} ^{\ast p}$. The possibility of $n(j) = p - 1$ is
ruled out in the same way, so (cc)$\to $(c), as claimed.
\par
(ccc)$\to $(cc) Suppose finally that $\rho _{j} = \varphi _{1}
(\eta _{j}) ^{p}\eta _{j} ^{-p}$, for some $\eta _{j} \in E _{j}
^{\ast }$. Then we have $(N _{1}N _{2}) \ast \eta _{j} ^{\prime }
= \rho _{j}$, for every $\eta _{j} ^{\prime } \in E _{j} ^{\ast }$
satisfying the equality $\varphi _{2} (\eta _{j} ^{\prime })\eta
_{j} ^{\prime -1} = \eta _{j} ^{p}N _{E} ^{E _{j}} (\eta _{j})
^{-1}$, so the proofs of (5.3) and Proposition 5.6 are complete.
\par
\vskip0.75truecm
\centerline {\bf 6. Intermediate norms in noncyclic abelian
extensions of degree $p ^{2}$}
\par
\medskip
The purpose of this Section is to prove Lemma 4.3. Let $E$ be a
field, $p$ a prime number, $F _{1}$ and $F _{2}$ different
extensions of $E$ in $E (p)$ of degree $p$, $\sigma $ an
$E$-automorphism of $F _{1}$ of order $p$, and $L = F _{1}F _{2}$.
By Lemma 5.1, then $L/E$ is a noncyclic Galois extension, $[L\colon
E] = p ^{2}$ and $I(L/E) = \{E _{1}, \dots , E _{p-1}, F _{1}, F
_{2}, E, L\}$. Note also that $E _{1} \cap F _{1} = E$ and $N _{E}
^{F _{1}} (\beta _{1}) = N ^{L} _{E _{1}} (\beta _{1})$, for every
$\beta _{1} \in F _{1}$. Considering now the cyclic $E$-algebra $A
_{\rho } = (F _{1}/E, \sigma , \rho )$, for an arbitrary $\rho \in
E ^{\ast }$, and observing that $A _{\rho } \otimes _{E} E _{1}$
is $E _{1}$-isomorphic to $(L/E _{1}, \bar \sigma , \rho )$, where
$\bar \sigma $ is the unique $E _{1}$-automorphism of $L$
extending $\sigma $, one obtains from (1.2), (1.4) (ii) and the
$p$-quasilocal property of $E$ that $\rho \in N(L/E _{1})$, i.e.
$E ^{\ast } \subseteq N(L/E _{1})$. This, combined with
Proposition 5.6, proves Lemma 4.3 in the case of $p = 2$. Suppose
further that $p > 2$ and put $E _{p+\mu } = F _{\mu}\colon \ \mu =
0, 1$. Then Lemmas 5.2 and 5.3 imply the following statements:
\par
\medskip
(6.1) (i) $E$ does not contain a primitive $p$-th root of unity if
and only if $(L ^{\ast p} \cap E) = E ^{\ast p}$.
\par
(ii) If $E$ contains a primitive $p$-th root of unity, then
conditions (a) and (b) of Lemma 5.2 (ii) can be
restated as follows:
\par
(a) $(L ^{\ast p} \cap N(E _{i}/E)) = (E _{i} ^{\ast p}
\cap E)$, $i = 1, \dots , p + 1$;
\par
(b) $(L ^{\ast p} \cap E) \subseteq N(E _{i}/E)$, $i = 1,
\dots , p + 1$.
\par
\medskip
Assume now that $c \in E ^{\ast }$ and $\xi _{1}$ is an element of
$L$ of norm $N _{E _{1}} ^{L} (\xi _{1}) = c$. The idea of our
proof is to establish consecutively the existence of elements $\xi
_{2}, \dots , \xi _{p-1}$ of $L$ satisfying the equalities $N _{E
_{i}} ^{L} (\xi _{k}) = c\colon \ i = 1, \dots , k$, for each
index $k$, and also to show that $\xi _{p-1}$ can be chosen so
 as to satisfy condition (ii) of Proposition 5.6. To implement
this we need additional information about the norms $N _{E _{j}}
^{L} (\xi _{k})\colon \ j = k + 1, \dots , p + 1$, for $k = 1, \dots ,
p - 1$. It is contained in the following four lemmas and seems to
be of independent interest.
\par
\medskip
{\bf Lemma 6.1.} {\it Let $E$ be a field , $p$ an odd prime number,
$L/E$ a noncyclic Galois extension of degree $p ^{2}$, $E _{1},
\dots , E _{p+1}$ the extensions of $E$ in $L$ of degree $p$,
$\varphi _{1}$ and $\varphi _{2}$ generators of ${\cal G}(L/E _{1})$
and ${\cal G}(L/E _{2})$, respectively, $N _{i} = \varphi _{i} - 1
\in \hbox{\Bbb Z} [{\cal G}(L/E)]\colon \ i =$ $1, 2$, and $\gamma $
an element of $L$ of norms $N _{E _{1}} ^{L} (\gamma ) = N _{E _{2}}
^{L} (\gamma ) = c$, for some $c \in E ^{\ast }$. Then $N _{E _{j}}
^{L} (\gamma ) = c.(N _{1}N _{2}) \ast \nu _{j}$, for some $\nu _{j}
\in E _{j} ^{\ast }$ and each index $j \ge 3$. Moreover:}
\par
(i) {\it If $p = 3$, then $c \in N(E _{3}/E)N(E _{4}/E)$;}
\par
(ii) {\it If $p \ge 5$, then $N _{E} ^{E _{3}} (\nu _{3})N _{E} ^{E
_{j}} (\nu _{j}) ^{-1}$ is contained in $L ^{\ast p}$, for $j = 3,
\dots , p + 1$; in addition, if $E$ is $p$-quasilocal, then there
exists $\gamma ^{\prime } \in L$, such that $N _{E _{i}} ^{L}
(\gamma ^{\prime }) = c\colon \ i = 1,$ $2, 3$.}
\par
\medskip
{\it Proof.} One can assume without loss of generality that, for
each index $j \ge 3$, ${\cal G}(L/E _{j})$ is generated by the element
$\varphi _{1}\tau _{j} ^{-1}$, where $\tau _{j} = \varphi _{2}
^{j-2}$. It is verified by direct calculations that the double
product $w _{j} (\lambda ) = \prod _{n=1} ^{p-1} (\prod _{u=1}
^{n} (\varphi _{1} ^{n-u} \tau _{j} ^{u}) (\lambda ))$ satisfies
the equality
$$(\varphi _{1}\tau _{j} ^{-1}) (w _{j} (\lambda ))w _{j}
(\lambda ) ^{-1} = \prod _{n=1} ^{p-1} (\varphi _{1} ^{n} (\lambda
)\tau _{j} ^{n} (\lambda ) ^{-1}) = N _{E _{1}} ^{L} (\lambda )N _{E
_{2}} ^{L} (\lambda ) ^{-1},$$ for any $\lambda \in L ^{\ast }$.
Similarly, one obtains that $\tau _{j} (w _{j} (\lambda )) = \prod
_{n=2} ^{p} (\prod _{u=2} ^{n} (\varphi _{1} ^{n-u}\tau _{j} ^{u})
(\lambda ))$ and
$$\tau _{j} (w _{j} (\lambda ))w _{j} (\lambda ) ^{-1} = [\prod
_{u=2} ^{p} (\varphi _{1} ^{p-u}\tau _{j} ^{u}) (\lambda )].[\prod
_{n=1} ^{p-1} (\varphi _{1} ^{n-1}\tau _{j}) (\lambda )] ^{-1} =$$
$$[\prod _{u=1} ^{p} (\varphi _{1} ^{p-u}\tau _{j} ^{u}) (\lambda
)].[\prod _{n=1} ^{p} (\varphi _{1} ^{n-1}\tau _{j}) (\lambda )]
^{-1} = N _{E _{j}} ^{L} (\lambda )N _{E _{1}} ^{L} (\tau _{j}
(\lambda )) ^{-1}.$$
As $L/E$ is abelian, this means that $\tau _{j}
(w _{j} (\lambda )) w _{j} (\lambda ) ^{-1} = N _{E _{j}} ^{L}
(\lambda )\tau _{j} (N _{E _{1}} ^{L} (\lambda ) ^{-1})$. These
results show that $\tau _{j} (w _{j} (\gamma ))w _{j} (\gamma )
^{-1} = N _{E _{j}} ^{L} (\gamma )N _{E _{1}} ^{L} (\gamma ) ^{-1}$
and $(\varphi _{1}\tau _{j} ^{-1}) (w _{j} (\gamma )) = w _{j}
(\gamma )$, i.e. $w _{j} (\gamma ) \in E _{j}$, for $j = 3, \dots
, p + 1$. Note also that
$$N _{E} ^{E _{j}}
(w _{j} (\gamma )) = N _{E _{1}} ^{L} (w _{j} (\gamma )) = \prod
_{n=1} ^{p-1} (\prod _{u=1} ^{n} N _{E _{1}} ^{L} ((\varphi _{1}
^{n-u}\tau _{j} ^{u}) (\gamma)) = c ^{p(p-1)/2}.$$
Observing now that $\varphi _{1}$ induces on $E _{j}$ an
$E$-automorphism of order $p$, and applying Hilbert's Theorem 90
as well as the inequality $p > 2$, one concludes that there exists
an element $\xi _{j} \in E _{j} ^{\ast }$ satisfying the
conditions $w _{j} (\gamma ) = c ^{(p-1)/2}\varphi _{1} (\xi
_{j})\xi _{j} ^{-1}$ and $N _{E _{j}} ^{L} (\gamma ) = c.(N
_{1}(\tau _{j} - 1)) \ast \xi _{j}$. Therefore, we have $N _{E
_{j}} ^{L} (\gamma ) = c.(N _{1}N _{2}) \ast \nu _{j}$,
where $\nu _{j} = \prod _{i=0} ^{j-3} \varphi _{2} ^{i} (\xi
_{j})$, for each index $j \ge
3$. Thus Lemma 6.1 (i) reduces to a special case of Proposition
5.6. Similarly, the conclusions of Lemma 6.1 (ii) are contained in
the following lemma.
\par
\medskip
{\bf Lemma 6.2.} {\it Assume that $E, p, L, E _{1}, \dots , E
_{p+1}$ are defined as in Lemma 6.1, $\varphi _{n}$ is a generator
of ${\cal G}(L/E _{n})$ and $N _{n} = \varphi _{n} - 1$, for every index
$n$. Let $k$ be an integer with $1 \le k < p - 1$, and suppose that
$\alpha \in L$ is of norms $N _{E _{i}} ^{L} (\alpha ) = c, \ i = 1,
\dots , k$, for a given $c \in E ^{\ast }$. Let also $N _{E _{j}}
^{L} (\alpha ) = c.(\prod _{i=1} ^{k} N _{i}) \ast \mu _{j}$, for
some $\mu _{j} \in E _{j} ^{\ast }$ and any $j \ge k + 1$. Then:}
\par
(i) {\it The products $N _{E} ^{E _{k+1}} (\mu _{k+1})N _{E} ^{E
_{j}} (\mu _{j}) ^{-1}\colon \ j = k + 1, \dots , p + 1$, are
contained in $L ^{\ast p}$; furthermore, if $k = 1$, then they lie
in $E ^{\ast p}$;}
\par
(ii) {\it If $E$ is $p$-quasilocal, then there exists an element
$\alpha ^{\prime } \in L$, such that $N _{E _{i'}} ^{L} (\alpha
^{\prime }) = c\colon $ $i ^{\prime } = 1, \dots , k + 1$.}
\par
\medskip
{\it Proof.} We begin with the latter assertion of Lemma 6.2 (i).
First we show that it suffices to consider the special case where
Br$(E) _{p} = \{0\}$. By (1.6) (ii), there exists an extension
$\Lambda $ of $E$, such that Br$(\Lambda ) _{p} = \{0\}$ and $E$ is
algebraically closed in $\Lambda $. Denote by $\widetilde E _{1},
\dots , \widetilde E _{p+1}$ and $\widetilde L$ the tensor products
$E _{1} \otimes _{E} \Lambda , \dots , E _{p+1} \otimes _{E} \Lambda
$ and $L \otimes _{E} \Lambda $, respectively. It is clear from
Galois theory and the equality $L _{\rm sep} \cap \Lambda = E$ that
$\Lambda ^{\ast p} \cap E = E ^{\ast p}$, $\widetilde L/\Lambda
$ is a noncyclic abelian extension of degree $p ^{2}$ and
$\widetilde E _{1}, \dots , \widetilde E _{p + 1}$ are the
extensions of $\Lambda $ in $\widetilde L$ of degree $p$. In
addition, it is easily verified that $N _{\Lambda } ^{\widetilde E
_{n}} (\eta _{n} \otimes _{E} 1) = N _{E} ^{E _{n}} (\eta _{n})
\otimes _{E} 1$ and $N _{\widetilde E _{n}} ^{\widetilde L} (\eta
\otimes _{E} 1) = N _{E _{n}} ^{L} (\eta ) \otimes _{E} 1\colon \ 1
\le n \le p + 1$, $\eta _{n} \in E _{n}$ and $\eta \in L$.
These observations lead to the desired reduction. By Lemma 4.2
(ii), then $N _{E _{2}} ^{L} (e _{2}) = \mu _{2}$, for
some $e _{2} \in L$. This implies that $N _{E _{u}} ^{L} (\alpha
^{\prime }) = c, \ u = 1, 2$, where $\alpha ^{\prime } = \alpha
\varphi _{1} (e _{2} ^{-1})e _{2}$. Applying Lemma 6.1, one
obtains that $\mu _{j}N _{E _{j}} ^{L} (e _{2}) ^{-1} = c
_{j}\varphi _{2} (\nu _{j})\nu _{j} ^{-1}$, for some $c _{j} \in E
^{\ast }$, $\nu _{j} \in E _{j} ^{\ast }$, and each index $j \ge
3$. Hence, $N _{E} ^{E _{j}} (\mu _{j})N _{E} ^{L} (e _{2}) ^{-1}
= c _{j} ^{p}$, which proves the latter part of Lemma 6.2 (i).
\par
For the proof of the former one, it is now sufficient to consider
the special case where $k \ge 2$. Fix an index $j \ge k + 2$,
denote by $M(j)$ and $K(j)$ the sets $\{k + 2, \dots , p + 1\}
\setminus \{j\}$ and $\{1, \dots , p + 1\} \setminus \{k + 1,
j\}$, respectively, and put $\alpha _{j} = (\prod _{m \in M(j)}
N _{m}) \ast \alpha $. It is easily verified that $N _{E _{j'}}
^{L} (\alpha _{j}) = 1$, for all $j ^{\prime } \in K(j)$. Taking
also into account that $N _{E} ^{L} (\alpha _{j}) = 1$, and for
every $n \in K(j)$, $\varphi _{n}$ induces on $E _{k+1}$ and $E
_{j}$ automorphisms of order $p$, one obtains by applying Lemma
5.4 that $\alpha _{j} ^{p} = (\prod _{n \in K(j)} N _{n}) \ast
(\mu _{k+1}\mu _{j})$. In view of (5.2), this yields $\alpha _{j}
^{p} = N _{E} ^{E _{k+1}} (\mu _{k+1}) ^{m(k+1)}N _{E} ^{E _{j}}
(\mu _{j}) ^{m(j)}\theta _{k+1} ^{p}\theta _{j} ^{p}$, for some
integers $m(k + 1)$ and $m(j)$ not divisible by $p$, and some
$\theta _{k+1} \in E _{k+1} ^{\ast }$, $\theta _{j} \in E _{j}
^{\ast }$. Therefore, the former statement of Lemma 6.2 (i) will
be proved, if we show that $p$ divides $m(k+1) + \ m(j)$. For each
index $n \ge 3$, denote by $l(n)$ the unique integer satisfying
the conditions $1 \le l(n) < p$ and $\varphi _{1}\varphi _{2}
^{-l(n)} \in {\cal G}(L/E _{n})$. Using the fact that $\varphi _{n} ^{l}
(\beta )\beta ^{-1} = \varphi _{n} (\beta _{l})\beta _{l} ^{-1}$,
where $\beta _{l} = \prod _{u=1} ^{l} \varphi _{n}  ^{u-1} (\beta
)$, for each $l \in \hbox{\Bbb N}$ and $\beta \in (E _{k+1} ^{\ast
} \cup E _{j} ^{\ast })$, one concludes that it suffices to
consider the special case where $\varphi _{n} = \varphi
_{1}\varphi _{2} ^{-l(n)}, n = 3, \dots , p + 1$. Then $\varphi
_{1} (\lambda _{\nu }) = \varphi _{2} ^{l(\nu )} (\lambda _{\nu
})$ and $\varphi _{n} (\lambda _{\nu }) = \varphi _{2} ^{l(\nu ) -
l(n)} (\lambda _{\nu })$, for $\nu \ge 3$ and $\lambda _{\nu } \in
E _{\nu }$. Hence, by Lemma 5.5 (i) and the inequality $k \ge 2$,
one can take as $m(k+1)$ and $m(j)$ the products
$$(6.2) \ l(k+1).1.\prod _{n \in L(j)} (l(k+1) - l(n)) \
{\rm and} \ l(j).1.\prod _{n \in L(j)} (l(j) - l(n)),$$
respectively, for a suitable choice of $\theta _{k+1}$ and $\theta
_{j}$, where $L(j) = \{3, \dots , p + 1\} \setminus $
\par \noindent
$\{k + 1, j\}$. Observe also that $\bar m(k+1) \equiv \bar m(j)
\equiv (p-1)! \ {\rm mod} \ p$, for $\bar m(k+1) =$ $(l(k+1) -
l(j))m(k+1)$ and $\bar m(j) = (l(j) - l(k+1))m(j)$. This implies
that
\par \noindent
$p \vert (m(k+1) + m(j))$ and so proves the
former assertion of Lemma 6.2 (i). The rest of our proof relies on
the following statements:
\par
\medskip
(6.3) (i) If $k = 1$, then $N _{E} ^{E _{2}} (\mu _{2}) \in N(E
_{3}/E)$;
\par
(ii) If $k \ge 2$ and $j \ge k + 2$, then there exist elements
$\lambda _{j} \in E _{k+1} ^{\ast }$ and $\omega _{j} \in E _{j}
^{\ast }$, such that $\lambda _{j} ^{p} \in E, \omega _{j} ^{p}
\in E$ and $N _{E} ^{E _{k+1}} (\lambda _{j}\mu _{k+1})N _{E} ^{E
_{j}} (\omega _{j}\mu _{j}) ^{-1} \in E ^{\ast p}$; in addition, $N
_{E _{k+1}} ^{L} (\alpha ) = c.(\prod _{i=1} ^{k} N _{i}) \ast
(\lambda _{j}\mu _{k+1})$ and $N _{E _{j}} ^{L} (\alpha ) =
c.(\prod _{i=1} ^{k} N _{i}) \ast (\omega _{j}\mu _{j})$.
\par
\medskip
Statement (6.3) (i) follows at once from the latter part of Lemma
6.2 (i), and by the second part of Lemma 5.5 (ii), $(N _{1}N _{2})
\ast t _{j'} = 1$ whenever $1 \le j ^{\prime } \le p + 1$, $t
_{j'} \in E _{j'} ^{\ast }$ and $t _{j'} ^{p} \in E ^{\ast }$. This
allows us to deduce (6.3) (ii) from Lemmas 5.2 and 5.3.
\par
We are now in a position to prove Lemma 6.2 (ii). Statement (6.3)
indicates that $N _{E _{k+1}} ^{L} (\alpha ) = c.(\prod _{i=1} ^{k}
N _{i}) \ast \mu _{k+1} ^{\prime }$ and $N _{E} ^{E _{k+1}} (\mu
_{k+1} ^{\prime }) \in N(E _{k+2}/E)$, for some $\mu _{k+1} ^{\prime
} \in $ $E _{k+1} ^{\ast }$. It is therefore clear from Lemma 4.2
(iii) that if $E$ is $p$-quasilocal, then $L$ contains an element $e
_{k+1}$ of norm $N _{E _{k+1}} ^{L} (e _{k+1}) = \mu _{k+1} ^{\prime
}$. In this case, the element $\alpha ^{\prime } = \alpha .(\prod
_{i=1} ^{k} N _{i}) \ast e _{k+1} ^{-1}$ satisfies the equations $N
_{E _{i'}} ^{L} (X) = c\colon \ i ^{\prime } = 1, \dots , k + 1$, so
Lemma 6.2 is proved.
\par
\medskip
{\bf Lemma 6.3.} {\it Assume that $E, p, L, E _{1}, \dots , E
_{p+1}$ are given as in Lemma 6.1, $c \in E ^{\ast }$, $k$ is an
integer with $1 \le k < p$, and $\ \alpha \in L$ satisfies the
equalities $N _{E _{n}} ^{L} (\alpha ) = c, n = 1, \dots , k$.
Suppose also that either $E$ does not contain a primitive $p$-th
root of unity or condition (6.1) (ii) (a) holds, and for each index
$n \le k$, let $\varphi _{n}$ be a generator of ${\cal G}(L/E _{n})$
and $N _{n} = \varphi _{n} - 1$. Then there exist elements $\mu
_{k+1} \in E _{k+1} ^{\ast }, \dots , \mu _{p+1} \in E _{p+1} ^{\ast
}$, such that
\par
\medskip \noindent
$(6.4) \ N _{E _{j}} ^{L} (\alpha ) = c.(\prod _{n=1} ^{k} N _{n})
\ast \mu _{j}$, $j = k + 1, \dots , p + 1$.
\par
\medskip \noindent
Moreover,}
\par
(i) {\it If $k < p - 1$, then the $(p + 1 - k)$-tuple $\bar \mu =
(\mu _{k+1}, \dots , \mu _{p+1})$ can be chosen so that $N _{E}
^{E _{k+1}} (\mu _{k+1})E ^{\ast p}$ $= N _{E} ^{E _{j}} (\mu
_{j})E ^{\ast p}$, for each index $j \ge k + 1$; in this case, the
co-set $N _{E} ^{E _{k+1}} (\mu _{k+1})E ^{\ast p}$ depends on
$\alpha $ and $\varphi _{1}, \dots , \varphi _{k}$ but not on the
choice of $\mu _{k+1}$;}
\par
(ii) {\it If $k < p - 1$ and $E$ is $p$-quasilocal, then $\mu
_{k+1}, \dots , \mu _{p+1}$ have the properties required by (i) if
and only if $\mu _{j} \in N(L/E _{j}), j = k + 1, \dots , p + 1$.}
\par
\medskip
{\it Proof.} First we prove the existence of elements $\mu _{j} \in
E _{j} ^{\ast }$, $j = k + 1, \dots , p + 1$, satisfying (6.4) and
with the properties required by the former statement of Lemma 6.3
(i). If $k = 1$, this is covered by Lemma 6.2, since then its second
hypothesis follows from Lemma 5.4. Henceforth, we consider the
special case of $k \ge 2$, assuming that our assertions are valid
for $k - 1$ and each element of $L$ of norm $c$ over $E _{n}\colon \
n = 1, \dots , k - 1$. This, applied to $\alpha $, enables one to
deduce from Lemma 5.5 (ii) the existence of elements $\tilde \mu
_{j} \in E _{j} ^{\ast }$, $j = k + 1, \dots , p + 1$, such that $N
_{E _{j}} ^{L} (\alpha ) = c.(\prod _{n=1} ^{k-1} N _{n}) \ast
\tilde \mu _{j}$ and $N _{E} ^{E _{j}} (\tilde \mu _{j}) \in E _{k}
^{p}$, for each index $j$. In view of (6.1) (i), (ii) (a) and Kummer
theory, this implies that $N _{E} ^{E _{j}} (\tilde \mu _{j}) \in E
^{\ast p}$. Hence, by Lemma 5.5 (iii), $N _{E _{j}} ^{L} (\alpha ) =
c.(\prod _{n=1} ^{k} N _{n}) \ast \mu _{j}$, for some $\mu _{j} \in
E _{j} ^{\ast }$ and every $j \ge k + 1$. Furthermore, it follows
from (6.1) (i) and Lemma 6.2 (i) that if $k < p - 1$ and $E$ does
not contain a primitive $p$-th root of unity, then $\mu _{k+1},
\dots , \mu _{p+1}$ have the properties required by the former
statement of Lemma 6.3 (i). Suppose now that $k < p - 1$ and (6.1)
(ii) (a) holds (with $E$ containing a primitive $p$-th root of
unity). Then one obtains from (6.3) and Lemma 6.2 that $\bar \mu $
can be fixed so that $N _{E} ^{E _{k+1}} (\mu _{k+1})N _{E} ^{E
_{k+2}} (\mu _{k+2}) ^{-1} \in E ^{\ast p}$ and $N _{E} ^{E _{k+1}}
(\mu _{k+1})N _{E} ^{E _{j}} (\mu _{j}) ^{-1} \in E _{k+1} ^{\ast
p}$, for $j = k + 3, \dots , p + 1$. We show that $N _{E} ^{E
_{k+1}} (\mu _{k+1})N _{E} ^{E _{j}} (\mu _{j}) ^{-1} \in E ^{\ast
p}, \ j \ge k + 3$. Statement (6.3) (ii) and our choice of $\bar \mu
$ guarantee that $N _{E} ^{E _{k+1}} (\mu _{k+1}) \in (N(E _{k+1}/E)
\cap $ $N(E _{k+2}/E))$ and for each $j \ge k + 3$, there exists
$\lambda _{j} \in E _{k+1} ^{\ast }$, such that
\par \noindent
$N _{E} ^{E _{k+1}} (\lambda _{j}\mu _{k+1})N _{E} ^{E _{j}} (\mu
_{j}) ^{-1}$ lies in $E ^{\ast p}$ and $\lambda _{j} ^{p} \in E$. As
$N _{E} ^{E _{k+1}} (\lambda _{j}) = \lambda _{j} ^{p}$, and by
Lemma 1.1, $N _{E} ^{E _{k+1}} (\mu _{k+1}) \in N(E _{j}/E)$, this
means that $\lambda _{j} ^{p} \in (L ^{\ast p} \cap N(E _{k+1}/E)
\cap $ $N(E _{j}/E))$. Moreover, by (6.1) (ii) (a), $\lambda _{j}
^{p} \in (E _{k+1} ^{\ast p} \cap E _{j} ^{\ast p})$. Since, by
Kummer theory, $E _{k+1} ^{\ast p} \cap E _{j} ^{\ast p} = E ^{\ast
p}$, this yields $\lambda _{j} \in E ^{\ast }$, for $j = k + 3,
\dots , p + 1$, which completes the proof of the existence part of
Lemma 6.3.
\par
Assume further that the elements $\mu _{k+1} \in E _{k+1}, \dots ,
\mu _{p+1} \in E _{p+1}$ satisfy (6.4) and have the properties
required by the former part of Lemma 6.3 (i), fix a $(p + 1 -
k)$-tuple $\bar u = (u _{k+1}, \dots , u _{p+1}) \in (E _{k+1}
\times \dots \times E _{p+1})$, and put $t _{j} = u _{j}\mu _{j}
^{-1}, j = k + 1, \dots ,$ $p + 1$. Clearly, we have $c.(\prod
_{n=1} ^{k} N _{n}) \ast u _{j} = N _{E _{j}} ^{L} (\alpha )$, for a
given index $j > k$, if and only if $(\prod _{n=1} ^{k} N _{n}) \ast
t _{j} = 1$. When $k = 1$ or $E$ does not contain a primitive $p$-th
root of unity, it follows from (6.1) (i), Lemma 5.5 (ii) and Galois
theory that this occurs if and only if $t _{j} \in E ^{\ast }$.
Therefore, in this case, the latter part of Lemma 6.3 (i) becomes
obvious, and Lemma 6.3 (ii) can be deduced from Lemma 4.2 (iii) (and
the inclusion $E ^{\ast } \subseteq N(L/E _{j})$ when $E$ is
$p$-quasilocal). Suppose now that $k \ge 2$ and (6.1) (ii) (a)
holds. Then it follows from Lemma 5.5 (ii) that $\bar u$ is a
solution to (6.4) if and only if $t _{j} ^{p} \in E ^{\ast }$, $j =
k + 1, \dots , p + 1$. When $\bar u$ is a solution, one sees that it
has the properties required by the former assertion of Lemma 6.3 (i)
if and only if all $t _{j}$ lie in $E$. This proves Lemma 6.3 (i).
At the same time, as above, one obtains from Lemma 4.2 (iii) that if
$E$ is $p$-quasilocal and $\bar u$ satisfies (6.4), then $u _{j'}
\in N(L/E _{j'})$, for a given index $j ^{\prime } \ge k + 1$, if
and only if $t _{j'} \in E$. This completes the proof of Lemma 6.3.
\par
\medskip
{\bf Lemma 6.4.} {\it Let $E$ be a field containing a primitive
$p$-th root of unity, and let $L/E$, $p$ and $E _{1}, \dots , E
_{p+1}$ satisfy (6.1) (ii) (b) and the conditions of Lemma 6.1.
Suppose that $k$ is an integer with $2 < k < p$, $\xi _{k}$ is a
$p$-th root in $E _{k}$ of an element $a _{k} \in (E ^{\ast }
\setminus E ^{\ast p})$, $\alpha \in L$ is of norms $N _{E _{i}}
^{L} (\alpha ) = c\colon \ i = 1, \dots , k$, for some $c \in E
^{\ast }$, and for each index $n$, $\varphi _{n}$ is a generator of
${\cal G}(L/E _{n})$ and $N _{n} = \varphi _{n} - 1$. Then there
exists an integer $\nu (k, \alpha )$ and a $(p + 1 - k)$-tuple $\bar
\mu = (\mu _{k+1}, \dots , \mu _{p+1}) \in (E _{k+1} \times \dots
\times $ $E _{p+1})$, such that $0 \le \nu(k, \alpha ) < p$, $N _{E
_{j}} ^{L} (\alpha ) = c.(\prod _{i=1} ^{k-1} N _{i}) \ast \mu _{j}$
and $N _{E} ^{E _{j}} (\mu _{j}) = a _{k} ^{\nu (k, \alpha )}$, for
$j = k + 1, \dots , p + 1$. Moreover,}
\par
(i) $\nu (k, \alpha )$ {\it does not depend on the choice of $\bar
\mu $ but is uniquely determined by $k, \alpha $ and $\varphi
_{1}, \dots , \varphi _{k-1}$;}
\par
(ii) $\nu (k, \alpha ) = 0$ {\it if and only if there are elements
$\lambda _{k+1} \in E _{k+1}, \dots , \lambda _{p+1} \in E _{p+1}$,
such that $N _{E _{j}} ^{L} (\alpha ) = c.(\prod _{i=1} ^{k} N
_{i}) \ast \lambda _{j}$, $j = k + 1, \dots , p + 1$;}
\par
(iii) {\it If $E$ is $p$-quasilocal, then $\xi _{k} \in N(L/E
_{k})$, and for each integer $m $ with $0 \le m < p$, there exists
$\alpha _{m} \in L$, such that $N _{E _{i}} ^{L} (\alpha _{m}) = c,
i = 1, \dots , k$, and $\nu (k, \alpha _{m})$ $= m$.}
\par
\medskip
{\it Proof.} Note first that it suffices to establish the
existence of an integer $\nu (k, \alpha )$ and of elements $\mu
_{k+1} \in E _{k+1}, \dots ,$ $\mu _{p+1} \in E _{p+1}$, such that
$0 \le \nu (k, \alpha ) < p$ and for every index $j \ge k + 1$, $N
_{E _{j}} ^{L} (\alpha ) = c.(\prod _{i=1} ^{k-1} N _{i}) \ast \mu
_{j}$ and $N _{E} ^{E _{j}} (\mu _{j}) = a _{k} ^{\nu (k, \alpha
)}$. Indeed, then Lemma 6.4 (i) can be deduced from Kummer theory
and Lemma 5.5 (ii), and Lemma 6.4 (ii) follows from Hilbert's
Theorem 90. When $E$ is $p$-quasilocal, Lemma 4.2 (iii) and the
inclusions $(L ^{\ast p} \cap E) \subseteq N(E _{n}/E)$ for $n = 1,
\dots , p + 1$, ensure the existence of an element $e _{k} \in L$ of
norm $N _{E _{k}} ^{L} (e _{k}) = \xi _{k}$. Using the inequality $k
\ge 3$, one easily verifies that
$$N _{E _{i}} ^{L} ((\prod _{u=1} ^{k-1} N _{u}) \ast e _{k} ^{m})
= (\prod _{u=1} ^{k-1} N _{u}) \ast N _{E _{i}} ^{L} (e _{k}) ^{m}
= 1\colon \ i = 1, \dots , k; \ m \in \hbox{\Bbb Z}.$$
This implies that the element $\alpha _{m} = \alpha .(\prod _{u=1}
^{k-1} N _{u}) \ast e _{k} ^{m-\nu (k,\alpha )} $ has the
properties required by Lemma 6.4 (iii), for $m = 0, 1, \dots , p - 1$.
\par
We turn to the main part of the proof of the lemma. Suppose first
that $k = 3$. By Lemma 6.1, then there exist $\nu _{4} \in E _{4}
^{\ast }, \dots , \nu _{p+1} \in E _{p+1} ^{\ast }$, such that $N
_{E _{j}} ^{L} (\alpha ) = c.(N _{1}N _{2})$ $\ast \nu _{j}\colon \
j = 4, \dots , p + 1$. Note also that by Lemma 6.2, $N _{E} ^{E
_{j}} (\nu _{j}) \in L ^{\ast p}$, for each $j \ge 4$. This enables
one to deduce from Lemmas 5.2 and 5.5 (ii) that $\nu _{j}$ can be
chosen so as to satisfy $N _{E} ^{E _{j}} (\nu _{j}) = a _{3}
^{n(j)}$, for some $n(j) \in \hbox{\Bbb Z}$ with $0 \le n(j) < p$.
We show that $n(j) = n(4)$ for $j = 5, \dots , p + 1$. As in the
proof of the latter assertion of Lemma 6.2 (i), our considerations
reduce to the special case in which Br$(E) _{p} = \{0\}$. By Lemma
4.2 (ii), then $E _{p+1}$ contains an element $\theta _{p+1}$ of
norm $N _{E} ^{E _{p+1}} (\theta _{p+1}) = c$. This implies that $N
_{E _{n}} ^{L} (\theta _{p+1}) = c\colon \  n = 1, \dots , p$, and
$N _{E _{i}} ^{L} (\alpha \theta _{p+1} ^{-1}) = 1\colon \ i = 1, 2,
3$. Hence, by Hilbert's Theorem 90, $\alpha \theta _{p+1} ^{-1} =
\varphi _{1} (\tilde \alpha )\tilde \alpha ^{-1}$, for some $\tilde
\alpha \in L ^{\ast }$. Moreover, it follows from Galois theory and
these facts that the norms $N _{E _{2}} ^{L} (\tilde \alpha ) :=
\rho _{2}$ and $N _{E _{3}} ^{L} (\tilde \alpha ) := \rho _{3}$ lie
in $E$. Since $N _{E} ^{L} (\tilde \alpha ) = N _{E} ^{E _{i}} (\rho
_{i}) = \rho _{i} ^{p}\colon \ i = 2, 3$, the element $\varepsilon =
\rho _{2}\rho _{3} ^{-1}$ is a $p$-th root of unity. Choose $\tilde
\theta _{p+1}$ from $E _{p+1}$ so that $N _{E} ^{E _{p+1}} (\tilde
\theta _{p+1}) = \rho _{2}$. Then $N _{E _{n}} ^{L} (\tilde \theta
_{p+1}) = \rho _{2}\colon \ n \le p$, and by Hilbert's Theorem 90,
the equation $\varphi _{2} (Y)Y ^{-1} =$ $\tilde \alpha \tilde
\theta _{p+1} ^{-1}$ has a solution $\bar \alpha \in L ^{\ast }$.
Clearly, we have $\varphi _{2} (N _{E _{3}} ^{L} (\bar \alpha ))N
_{E _{3}} ^{L} (\bar \alpha ) ^{-1} = \varepsilon ^{-1}$, which
means that $N _{E _{3}} ^{L} (\bar \alpha ) = \omega \xi _{3} ^{\mu
}$, for some $\omega \in E ^{\ast }$, $\mu \in \{0, 1, \dots , p -
1\}$. By Lemma 4.2 (ii), there exists $\lambda \in L ^{\ast }$ of
norm $N _{E _{3}} ^{L} (\lambda ) = \xi _{3} ^{\mu }$. This implies
that $N _{E _{3}} ^{L} (\bar \alpha \lambda ^{-1}) = \omega =$
$\varphi _{2} (N _{E _{3}} ^{L} (\bar \alpha \lambda ^{-1}))$. As
$\tilde \alpha = \tilde \theta _{p+1}\varphi _{2} (\bar \alpha )\bar
\alpha ^{-1}$, one also sees that $N _{E _{2}} ^{L} (\tilde \alpha
_{\lambda }) = N _{E _{3}} ^{L} (\tilde \alpha _{\lambda }) = \rho
_{2}$, where $\tilde \alpha _{\lambda } = \tilde \alpha \varphi _{2}
(\lambda ^{-1})\lambda $. Therefore, by Lemma 6.1, there are $c _{j}
^{\prime } \in E _{j}$ such that
$$N _{E _{j}} ^{L} (\tilde \alpha _{\lambda }) = \rho _{2}.(N _{2}N
_{3}) \ast c _{j} ^{\prime }, \ {\rm for} \ 1 \le j \le p + 1, \ j
\neq 2, 3.$$
Observing also that the equation $(\prod _{u=1} ^{p-1} N _{u})
\ast X _{p+1} = \theta _{p+1} ^{p}c ^{-1}$ has a root in $E
_{p+1}$ (see the proof of implication (i)$\to $(ii) of Proposition
5.6), one concludes that the elements $\alpha .(N _{1}N _{2}) \ast
\lambda ^{-1} = \theta _{p+1}\varphi _{1} (\tilde \alpha _{\lambda
})\tilde \alpha _{\lambda } ^{-1}$ and $c$ satisfy the conditions
of Lemma 6.2, for $k = 3$. In other words, $N _{E _{i}} ^{L}
(\alpha .(N _{1}N _{2}) \ast \lambda ^{-1}) = c\colon \ i = 1, 2,
3$, and
$$N _{E _{j}} ^{L} (\alpha .(N _{1}N _{2}) \ast \lambda ^{-1}) =
c.(N _{1}N _{2}) \ast (\nu _{j}N _{E _{j}} ^{L} (\lambda ^{-1})) =
c.(N _{1}N _{2}N _{3}) \ast c _{j}, \ j = 4, \dots , p + 1,$$ where
$c _{j} = c _{j} ^{\prime }\colon \ j \le p$ and $c _{p+1}$ is some
element of $E _{p+1} ^{\ast }$. Hence, by Lemma 5.5 (ii), one can
find an element $\delta _{j} \in E _{j} ^{\ast }$, such that $\delta
_{j} ^{p} \in E$ and $\nu _{j}N _{E _{j}} ^{L} (\lambda ^{-1}) =
\delta _{j}\varphi _{3} (c _{j})c _{j} ^{-1}$. The obtained result
indicates that $N _{E} ^{E _{j}} (\nu _{j})N _{E} ^{L} (\lambda
^{-1}) = a _{3} ^{n(j)-\mu } = \delta _{j} ^{p}$, i.e. $a _{3}
^{n(j)-\mu } \in $ $E _{j} ^{\ast p}$. In view of Kummer theory and
the assumptions on $a _{3}$, $\mu $ and $n(j)$, this means that
$n(j) = \mu $, for $j = 4, \dots , p + 1$, as claimed.
\par
Assume now that $k > 3$ and the conclusions of the lemma are valid
for $k - 1$, every subset $\{\Phi _{1}, \dots , \Phi _{k-1}\}$ of
$\{E _{1}, \dots E _{k}\}$, and each pair $(\alpha ^{\prime }, c
^{\prime }) \in (L ^{\ast } \times E ^{\ast })$ satisfying the
equalities $N _{\Phi _{i}} ^{L} (\alpha ^{\prime }) = c ^{\prime },
i = 1, \dots , k - 1$. This, applied to $\alpha $ and $(E _{1},
\dots , E _{k-1})$, implies the existence of elements $\mu _{j}
^{\prime } \in E _{j} ^{\ast }\colon \ j = k + 1, \dots , p + 1$,
such that $c.(\prod _{i=1} ^{k-2} N _{i}) \ast \mu _{j} ^{\prime } =
N _{E _{j}} ^{L} (\alpha )$, $N _{E} ^{E _{j}} (\mu _{j} ^{\prime })
\in E _{k-1} ^{\ast p}$ and $N _{E} ^{E _{j}} (\mu _{j} ^{\prime })
= N _{E} ^{E _{k}} (\mu _{k} ^{\prime })$, for each index $j$ and
some $\mu _{k} ^{\prime } \in E _{k} ^{\ast }$ satisfying the
equality $(\prod _{i=1} ^{k-2} N _{i}) \ast \mu _{k} ^{\prime } =
1$. In view of Kummer theory and Lemma 5.5 (ii), this yields $N _{E}
^{E _{k}} (\mu _{k} ^{\prime }) \in E ^{\ast p}$, so it follows from
Lemma 5.5 (iii) that $N _{E _{j}} ^{L} (\alpha ) = c.(\prod _{i=1}
^{k-1} N _{i}) \ast \mu _{j}$, for some $\mu _{j} \in E _{j} ^{\ast
}$. By Lemma 6.2 and the equality $N _{E _{k}} ^{L} (\alpha ) = c$,
this means that $N _{E} ^{E _{j}} (\mu _{j}) \in L ^{\ast p}$.
Furthermore, it becomes clear from Lemma 5.2 that $\mu _{j}$ can be
chosen so that $N _{E} ^{E _{j}} (\mu _{j}) = a _{k} ^{n(j)}$, for
some $n(j) \in \hbox{\Bbb Z}$ with $0 \le n(j) \le p - 1$. It
remains to be seen that $n(j) = n(k+1)$, $j = k + 2, \dots , p + 1$.
As in the case of $k = 3$, we obtain that one may assume in addition
that Br$(E) _{p} = \{0\}$. By Lemma 4.2 (ii), then $E _{p+1}$
contains an element $\theta _{p+1}$ of norm $N _{E} ^{E _{p+1}}
(\theta _{p+1}) = c$. Observing that $N _{E _{n}} ^{L} (\theta
_{p+1}) = c$, $n = 1, \dots , p$, and $N _{E _{i}} ^{L} (\alpha
\theta _{p+1} ^{-1}) = 1$, $i = 1, \dots , k$, one deduces from
Hilbert's Theorem 90 that $\alpha \theta _{p+1} ^{-1} = \varphi _{1}
(\tilde \alpha )\tilde \alpha ^{-1}$, for some $\tilde \alpha \in L
^{\ast }$. Also, it follows from Galois theory that $N _{E _{i}}
^{L} (\tilde \alpha ) := \rho _{i}$ lies in $E ^{\ast }$, for $i =
2, \dots , k$. We show that the $\rho _{i}$'s are equal. Our
argument relies upon the fact that $N _{E} ^{L} (\tilde \alpha ) = N
_{E} ^{E _{i}} (\rho _{i}) = \rho _{i} ^{p}$, i.e. the elements
$\varepsilon _{i} = \rho _{i}\rho _{k} ^{-1}$ are $p$-th roots of
unity. By Lemma 4.2 (ii), there exists an element $\tilde \theta
_{p+1} \in E _{p+1}$, such that $N _{E} ^{E _{p+1}} (\tilde \theta
_{p+1}) = \rho _{k}$. Note also that $N _{E _{n}} ^{L} (\tilde
\theta _{p+1}) = \rho _{k}\colon \ n = 1, \dots , p$, $N _{E _{i}}
^{L} (\tilde \alpha \tilde \theta _{p+1} ^{-1}) = \varepsilon _{i}$,
$i = 2, \dots , k - 1$, and $N _{E _{k}} ^{L} (\tilde \alpha \tilde
\theta _{p+1} ^{-1}) = 1$. Hence, by Hilbert's Theorem 90, the
equation $\varphi _{k} (Z)Z ^{-1} =$ $\tilde \alpha \tilde \theta
_{p+1} ^{-1}$ has a solution $\bar \alpha \in L ^{\ast }$. Observing
that the norm $N _{E _{i}} ^{L} (\bar \alpha ) := \omega _{i}$
satisfies the equality $\varphi _{k} (\omega _{i})\omega _{i} ^{-1}
= \varepsilon _{i}$, one obtains from Lemma 5.5 (ii) that $N _{E}
^{L} (\bar \alpha ) = N _{E} ^{E _{i}} (\omega _{i}) =$ $\omega _{i}
^{p}$. As $k > 3$, Kummer theory and these calculations show that $N
_{E} ^{L} (\bar \alpha ) \in E ^{\ast p}$, which yields
consecutively $\omega _{i} \in E ^{\ast }$, $\varepsilon _{i} = 1$
and $\rho _{i} = \rho _{k}$, for all $i = 2, \dots , k - 1$.
Therefore, by hypothesis, one can find an integer $\nu (k, \tilde
\alpha )$ and elements $\delta _{1} ^{\prime } \in E _{1}$, $\delta
_{k+1} ^{\prime } \in E _{k+1}, \dots , \delta _{p+1} ^{\prime } \in
E _{p+1}$, such that
$$0 \le \nu (k, \tilde \alpha ) \le p - 1, N _{E _{j}} ^{L}
(\tilde \alpha ) = \rho _{k}(\prod _{i=2} ^{k-1} N _{i}) \ast
\delta _{j} ^{\prime } \ {\rm and} \ N _{E} ^{E _{j}} (\delta _{j}
^{\prime }) = a _{k} ^{\nu (k, \tilde \alpha )},$$
for $j = 1$ and $j \ge k + 1$. It has also been pointed out that,
by the proof of Proposition 5.6, the choice of $\theta _{p+1}$
ensures the solvability of the equation $(\prod _{u=1} ^{p-1} N
_{u}) \ast X _{p+1} = \theta _{p+1} ^{p}c ^{-1}$ over $E _{p+1}$.
In view of the equality $\alpha = \theta _{p+1}\varphi _{1}
(\tilde \alpha )\tilde \alpha ^{-1}$, these results yield
$$N _{E _{j}} ^{L} (\alpha ) = c.(\prod _{i=1} ^{k-1} N _{i}) \ast
\delta _{j} ^{\prime \prime }\colon \ j = k + 1, \dots , p + 1,$$
where $\delta _{j} ^{\prime \prime } = \delta _{j} ^{\prime }\colon
j \le p$, and $\delta _{p+1} ^{\prime \prime } = \delta _{p+1}
^{\prime }.(\prod _{z=k} ^{p-1} N _{z}) \ast \theta _{p+1} ^{\prime
}$, for any $\theta _{p+1} ^{\prime } \in E _{p+1}$ satisfying the
equality $(\prod _{u=1} ^{p-1} N _{u}) \ast \theta _{p+1} ^{\prime }
= \theta _{p+1} ^{p}c ^{-1}$. This implies that $(\prod _{i=1}
^{k-1} N _{i}) \ast $ $(\mu _{j}\delta _{j} ^{\prime \prime -1}) =
1$ and $N _{E} ^{E _{p+1}} (\delta _{p+1} ^{\prime \prime }) = N
_{E} ^{E _{p+1}} (\delta _{p+1} ^{\prime })$, so it follows from
Lemma 5.5 (ii) that $E _{j} ^{\ast p}$ contains the norm $N _{E} ^{E
_{j}} (\mu _{j}\delta _{j} ^{\prime \prime -1}) = a _{k} ^{n(j) -
\nu (k, \tilde \alpha )}$, for every $j > k$. It is now clear from
Kummer theory and the condition on $a _{k}$ that $p \vert (n(j) -
\nu (k, \tilde \alpha ))$. The obtained result and the assumptions
on $n(j)$ and $\nu (k, \tilde \alpha )$ indicate that $n(j) = \nu
(k, \tilde \alpha )$, $j = k + 1, \dots , p + 1$ (i.e. one may put
$\nu (k, \alpha ) = \nu (k, \tilde \alpha )$), which completes the
proof of Lemma 6.4.
\par
\medskip
{\it Proof of Lemma 4.3.} Assume that $p > 2$ and $F _{1}F _{2} =
L$, take $E _{1}, \dots , E _{p-1}$ as at the beginning of this
Section, put $E _{p+\mu } = F _{\mu }\colon \ \mu = 0, 1$, and fix
an arbitrary element $c \in E ^{\ast }$. We have already proved
that $N _{E _{1}} ^{L} (\xi _{1}) = c$, for some $\xi _{1} \in L$.
Hence, by Lemmas 5.4 and 6.2 (ii), $L$ contains an element
$\xi _{2}$ such that $N _{E _{i}} ^{L} (\xi _{2}) = c$, $i = 1,
2$. In view of Lemma 6.1, this proves Lemma 4.3 in the case of $p
= 3$. Suppose now that $p \ge 5$. By Lemma 6.1 (ii), then there
exists $\xi _{3} \in L$ of norms $N _{E _{i}} ^{L} (\xi _{3}) =
c$, $i = 1, 2, 3$. Combining finally Lemma 6.2 with Lemma 6.3 or
6.4 (and applying (6.1)), one obtains
consecutively the existence of elements $\xi _{4}, \dots , \xi
_{p-1} \in L$ such that $N _{E _{i}} ^{L} (\xi _{j}) = c$, $i = 1,
\dots j$, for each index $j $. Furthermore, one concludes that
$\xi _{p-1}$ can be chosen so as to satisfy condition (ii) of
Proposition 5.6. This shows that $c \in N(F _{1}/E)N(F _{2}/E)$,
so Lemma 4.3 is proved.
\par
\vskip.75truecm \centerline {\bf 7. Proof of Theorem 4.1}
\par
\medskip
First we complete the technical preparation for the proof of our
main result by showing that the class of $p$-quasilocal fields is
closed under the formation of cyclic extensions of degree $p$.
This is carried out in two steps stated as lemmas.
\par
\medskip
{\bf Lemma 7.1.} {\it Assume that $E$ is a $p$-quasilocal field
for a given prime number $p$, $L/E$ is a Galois extension of
degree $p ^{2}$, $F$ is an extension of $E$ in $L$ of degree $p$,
$\Delta \in d(F)$ and {\rm ind}$(\Delta ) = p$. Then $L$ is
embeddable in $\Delta $ as an $F$-subalgebra.}
\par
\medskip
{\it Proof.} As $F/E$ is cyclic, Lemma 4.2 (i) and Theorem 3.1
imply that $\Delta $ is similar over $F$ to $D \otimes _{E} {F}$,
for some $D \in d(E)$ of exponent $p ^{2}$. Since the $L$-algebras
$D \otimes _{E} {L}$ and $(D \otimes _{E} {F}) \otimes _{F} {L}$
are isomorphic (cf. [P, Sect. 9.4, Corollary a]), hence similar to
$\Delta \otimes _{F} {L}$, this means that the conclusion of Lemma
7.1 can be restated by saying that $L$ embeds in $D$ as an
$E$-subalgebra. In particular, by Theorem 3.1 (iii), it holds in
the case where $L$ is cyclic over $E$. Suppose further that $L/E$
is noncyclic, i.e. $L = MF$, where $M$ is a cyclic extension of
$E$ in $L$ of degree $p$, $M \neq F$. Also, let $\sigma $ be an
$E$-automorphism of $M$ of order $p$, $\tilde \sigma $ the unique
$F$-automorphism of $L$ extending $\sigma $, and $D _{1}$ the
underlying division algebra of the $p$-th tensor power of $D$
over $E$. Observing that $D _{1} \in d(E)$ and ind$(D _{1}) = p$,
one obtains from (1.4) and Lemma 4.3 that $D _{1} \cong (M/E,
\sigma , c)$, for some $c \in N(F/E)$. Let $\gamma $ be an
element of $F$ of norm $N _{E} ^{F} (\gamma ) = c$. As noted in
the proof of Lemma 4.2 (iii), then cor$_{F/E}$ maps $[(L/F,
\tilde \sigma , \gamma )]$ into $[(M/E, \sigma , c)]$. At the
same time, by [T, Theorem 2.5], we have cor$_{F/E} ([\Delta ]) =
[D _{1}]$. Hence, by the injectivity of cor$_{F/E}$ and the
equality $[\Delta \colon F]$ $= [(L/F, \tilde \sigma , \gamma
)\colon F] = p ^{2}$, $\Delta $ and $(L/F, \tilde \sigma , \gamma
)$ are $F$-isomorphic, which proves Lemma 7.1.
\par
\medskip
The application of the corestriction mapping in the proofs of Lemmas
4.2 (iii) and 7.1 was suggested by the referee (for somewhat longer
proofs relying only on general properties of crossed products, see
e.g., the cross-reference in the proof of [Ch6, (4.1) (iii)]).
\par
\medskip
{\bf Lemma 7.2.} {\it Let $E$ be a $p$-quasilocal field and $F$ a
cyclic extension of $E$ of degree $p$. Then $F$ is
$p$-quasilocal.}
\par
\medskip
{\it Proof.} For each $\chi \in X _{p} (F)$, denote by $L _{\chi
}$ the cyclic extension of $F$ in $F _{\rm sep}$ fixed by ${\rm
Ker}(\chi )$, and by $s$ the pairing $X _{p} (F) \times F ^{\ast }
\to \ _{p} {\rm Br}(F)$ defined as in the proof of Lemma 1.1.
Suppose also that $\tau \in {\cal G}(F/E)$ and take any automorphism
$\rho $ of $F _{\rm sep}$ extending $\tau $. Then $\tau $ acts on
$X _{p}(F)$ by $\tau (\chi ) (g) = \chi (\rho ^{-1}g\rho )$, for
all $g \in {\cal G} _{F}$. Because $X _{p} (F)$ is an abelian group of
exponent $p$, this allows us to view it as a module over the group
ring $\hbox{\Bbb F} _{p} [{\cal G}(F/E)]$. Note that for any $\chi \in X
_{p} (F)$ and $b \in F ^{\ast }$, we have $\tau (s(\chi , b)) =
s(\tau (\chi ), \tau (b))$. Observe that
$$(7.1) \ {\rm if} \ \tau (s(\chi , b)) = s(\chi , b){\rm ,} \
{\rm then} \ s(\chi - \tau (\chi ), \tau (b)) = s(\chi , \tau (b)b
^{-1}).$$
$${\rm For} \ s(\chi - \tau (\chi ), \tau (b)) = s(\chi , \tau
(b)) - s(\tau (\chi ), \tau (b)) = s(\chi , \tau (b)) - s(\chi ,
b) = s(\chi , \tau (b)b ^{-1}),$$
by the $\hbox{\Bbb Z}$-bilinearity of $s$. Now to prove Lemma 7.2
note that ${\cal G}(F/E)$ acts trivially on $_{p} {\rm Br} (F)$, by Lemma
4.2 (i) and (1.5). Consider a cyclic extension $L$ of $F$ in $F
_{\rm sep}$ of degree $p$. For the proof of the lemma, it suffices
to show that Br$(L/F) = \ _{p} {\rm Br}(F)$ (see (1.1) (i) and
(1.2) (ii)). The given field $L$ is $L _{\chi }$ for some $\chi
\in X _{p} (F)$. Define inductively $\chi _{1} = \chi ,$ $\chi
_{2} = \chi _{1} - \psi (\chi _{1}), \dots , \chi _{i+1} = \chi
_{i} - \psi (\chi _{i}), \dots $, where $\psi $ is a fixed
generator of ${\cal G}(F/E)$. As $X _{p} (F)$ is an $\hbox{\Bbb F} _{p}
[{\cal G}(F/E)]$-module and (by Lemma 1.2) $1 - \psi $ is nilpotent in
$\hbox{\Bbb F} _{p} [{\cal G}(F/E)]$, we have $\chi _{l} = 0$ for $l$
sufficiently large. Choose $k$ so that $\chi _{k+1} = 0$ but $\chi
_{k} \neq 0$. Since $\psi (\chi _{k}) = \chi _{k}$, $L _{\chi
_{k}}$ is Galois over $E$ (of degree $p ^{2}$). Hence, Br$(L
_{\chi _{k}}/F) = \ _{p}$Br$(F)$ by Lemma 7.1. But because
${\cal G}(F/E)$ acts trivially on $_{p} {\rm Br}(F)$, statement (7.1)
(with $\psi $ for $\tau $) shows that Br$(L _{\chi _{i}}/F)
\subseteq {\rm Br}(L _{\chi _{i-1}}/F)$, for each $i \ge 2$. In view of
the inclusion Br$(L/F) \subseteq \ _{p} {\rm Br}(F)$, this proves
that Br$(L/F) = {\rm Br}(L _{\chi _{k}}/F) = \ _{p} {\rm Br}(F)$,
as desired.
\par
\medskip
It is now easy to prove Theorem 4.1. Suppose first that $R$ is an
extension of $E$ in $E (p)$ of degree $p ^{n}$, for some $n \in
\hbox{\Bbb N}$, and fix an extension $U$ of $E$ in $R$ of degree
$p$. Clearly, $[R\colon U] = p ^{n-1}$, and by Lemma 7.2, $U$ is a
$p$-quasilocal field. This, combined with Lemma 4.2 (i) and the
equality $\pi _{R/E} = \pi _{R/U} \circ \pi _{U/E}$, enables one
to prove by induction on $n$ that $R$ has the properties required
by Theorem 4.1 (i)-(ii). Assuming that Br$(E) _{p} \neq \{0\}$,
fix an algebra $D \in d(E)$ of index divisible by $p$, and put
g.c.d.$([R\colon E], {\rm ind}(D)) = p ^{k}$. By Galois theory and
the subnormality of proper subgroups of finite $p$-groups, $R/E$
possesses an intermediate field $U _{k}$ such that $[U
_{k}\colon E] = p ^{k}$ and $U \subseteq U _{k}$. We show that
$U _{k}$ embeds in $D$ as an $E$-subalgebra. Let $D _{1}$ be the
underlying division algebra of $D \otimes _{E} U$. By Theorem 3.1
(i)-(iv), applied to $D/E$ and $D _{1}/U$, we have exp$(D) = {\rm
ind}(D)$ and ind$(D _{1}) = {\rm exp}(D _{1}) = {\rm exp}(D)/p =
{\rm ind}(D)/p$. At the same time, the equality $\pi _{U _{k}/E} =
\pi _{U _{k}/U} \circ \pi _{U/E}$ ensures that $D \otimes _{E} U
_{k}$ and $D _{1} \otimes _{U} U _{k}$ are similar over $U _{k}$.
Since $[U _{k}\colon U] = p ^{k-1}$, these observations and (1.2)
(ii) indicate that $U _{k}$ embeds in $D$ over $E$ if and only if
it embeds in $D _{1}$ over $U$. Now the embeddability of $U _{k}$
in $D$ is easily proved by induction on $k$. As $\pi _{R/E} = \pi
_{R/U _{k}} \circ \pi _{U _{k}/E}$, this result, statement (1.2) (ii)
and Theorem 4.1 (iv) imply Theorem 4.1 (iii).
\par
It remains to be seen that if $R$ is an infinite extension of $E$
in $E (p)$, then Br$(R) _{p} = \{0\}$. Let $\Delta \in d(R)$ be of
$p$-primary dimension. By (1.3), there exists an $R$-isomorphism
$\Delta \cong \Delta _{0} \otimes _{R _{0}} R$, for some finite
extension $R _{0}$ of $E$ in $R$, and some $\Delta _{0} \in d(R
_{0})$. The $R _{0}$-algebra $\Delta _{0}$ is split by $R$, since
$R _{0}$ is $p$-quasilocal, $\pi _{R _{0}/R} = \pi _{R _{0}'/R}
\circ \pi _{R _{0}/R _{0}'}$ for every intermediate field $R _{0}
^{\prime }$ of $R/R _{0}$, and since for each $m \in \hbox{\Bbb N}$,
$R$ contains as a subfield an extension $R _{m}$ of $R _{0}$ of
degree $p ^{m}$. Therefore, $\Delta = R$ and Br$(R) _{p} = \{0\}$,
so Theorem 4.1 is proved.
\par
\vskip0.75truecm \centerline{\bf 8. On the absolute Galois groups of
absolutely stable and of quasilocal fields}
\par
\medskip
Now we turn our attention to the residue fields of Henselian valued
absolutely stable fields with totally indivisible value groups.
Proposition 2.3 and [Ch1, Corollary 4.6] indicate that a perfect
field $E$ is isomorphic to such a residue field if and only if $E$
is quasilocal. Our next result characterizes nonreal perfect
quasilocal fields and almost perfect absolutely stable fields by
cohomological properties of the Sylow subgroups of their absolute
Galois groups. Supplemented in [Ch3, Sect. 3] by a similar treatment
of the formally real case, it shows that quasilocal fields form one
of the basic classes of absolutely stable fields.
\par
\medskip
{\bf Theorem 8.1.} {\it Let $E$ be a field, $\Pi (E)$ the set of
all prime numbers $p$ for which ${\rm cd} _{p} ({\cal G} _{E}) \neq 0$,
and $\{G _{p}\}$ a set of Sylow pro-$p$-subgroups of ${\cal G} _{E}$,
indexed by  $\Pi (E)$. Then:}
\par
(i) {\it If $E$ is quasilocal and nonreal, then $G _{p}$ is a
$p$-group of Demushkin type, for each $p \in \Pi (E)$; conversely,
if $E$ is perfect with $G _{p}$ a $p$-group of Demushkin type, for
every $p \in \Pi (E)$, then $E$ is nonreal and quasilocal.}
\par
(ii) {\it If $E$ is absolutely stable, then the cup-product mapping
of $H ^{1} (G _{p} ^{\prime }, \hbox{\Bbb F} _{p}) \times H ^{1} (G
_{p} ^{\prime }, \hbox{\Bbb F} _{p})$ into $H ^{2} (G _{p} ^{\prime
}, \hbox{\Bbb F} _{p})$ is surjective, for every open subgroup $G
_{p} ^{\prime }$ of $G _{p}$ and each $p \in \Pi (E)$; the converse
is true, provided that $E$ is almost perfect.}
\par
\medskip
The proof and the applications of Theorem 8.1 rely on the fact that
if $\widetilde E/E$ is a purely inseparable field extension, then $E
_{\rm sep} \otimes _{E} \widetilde E$ is a separable closure of
$\widetilde E$ and there exist group isomorphisms ${\cal G}
_{\widetilde E} \cong {\cal G} _{E}$ and ${\cal G}(\widetilde E
(p)/\widetilde E) \cong {\cal G}(E (p)/E)$, $p \in {\cal P}(E)$. For
instance, when $\widetilde E$ is a perfect closure of $E$ and $E$ is
taken as required by (4.2), it is thereby proved that $\widetilde E$
is quasilocal. Since ${\cal G}(F (p)/F)$ is a $p$-group of Demushkin
type, for every quasilocal nonreal field $F$ and each $p \in {\cal
P}(F)$ (see [Ch8, Sect. 3]), this allows us to view the first part
of (4.2) as a description of the spectrum of values of the main
cohomological invariants of quasilocal fields. Note that this
spectrum is much wider than the one in the case of local fields, and
more generally, of Henselian valued quasilocal fields with totally
indivisible value groups (see [Se1, Ch. II, 2.2 and 5.6], [Wa2,
Lemma 7] and the comment on Proposition 8.9).
\par
\medskip
{\it Proof of Theorem 8.1.} Our argument relies on the following lemma.
\par
\medskip
{\bf Lemma 8.2.} {\it The classes of absolutely stable fields and
of quasilocal fields are closed under the formation of algebraic
extensions.}
\par
\medskip
{\it Proof.} The assertion about the class of absolutely stable
fields follows at once from (1.3) (i) and (iii). Let now $E/E _{0}$
be an algebraic field extension, $W \in d(E)$, and let $F/E$ be a
cyclic extension of degree $n$ dividing ind$(W)$. By (1.3)
(i)-(ii), there is a finite extension $F ^{\prime }$ of $E _{0}$
in $F$ and a central division algebra $W ^{\prime }$ over the
field $E \cap F ^{\prime } := E ^{\prime }$, such that $F ^{\prime
}/E ^{\prime }$ is cyclic of degree $n$ and the $E$-algebras $W
^{\prime } \otimes _{E'}E$ and $F ^{\prime } \otimes _{E'}E$ are
isomorphic to $W$ and $F$, respectively. It is therefore clear
that if $E _{0}$ is quasilocal, then $F$ embeds in $W $ as an
$E$-subalgebra, which proves our assertion about the class of
quasilocal fields.
\par
\medskip
Statement (1.8), [Ch3, Proposition 3.1] and our next lemma
indicate that it is sufficient to prove Theorem 8.1 in the special
case where ${\cal G} _{E}$ is a pro-$p$-group, for some $p \in {\cal
P}(E)$.
\par
\medskip
{\bf Lemma 8.3.} {\it Let $E$ be a field, ${\cal P}$ the set
of prime numbers, and for each $p \in {\cal P}$, let $G _{p}$
be a Sylow pro-$p$-subgroup of ${\cal G} _{E}$ and $E _{p} = \{\alpha \in
E _{\rm sep}\colon \ \sigma _{p} (\alpha ) = \alpha $, $\sigma
_{p} \in G _{p}\}$. Then:}
\par
(i) $E$ {\it is absolutely stable if and only if $E _{p}$ have the
same property, for all $p \in {\cal P}$;}
\par
(ii) $E$ {\it is quasilocal if and only if so are $E _{p}\colon \
p \in {\cal P}$; this is the case if and only if $E _{p}$ is
$p$-quasilocal, for every $p \in {\cal P}$.}
\par
\medskip
{\it Proof.} Note first that Br$(E _{p}) = {\rm Br}(E _{p}) _{p} =
{\rm Br}(E _{\rm sep}/E _{p})$, for every $p \in {\cal P}$ (see
P, Sect. 13.5]). This, combined with Theorem 4.1 and Proposition
4.4, reduces the latter conclusion of Lemma 8.3 (ii) to a
consequence of the former one. We prove Lemma 8.3 (i) and the
former part of Lemma 8.3 (ii). Let $\overline E$ be an algebraic
closure of $E _{\rm sep}$, $F$ a finite extension of $E$ in
$\overline E$, and $F _{0} = F \cap E _{\rm sep}$. Consider an
algebra $D \in d(F)$ of $p$-power index, for some $p \in {\cal
P}$, and a cyclic extension $L/F$ of degree dividing ind$(D)$. It
follows from Sylow's theorem (cf. [Se1, Ch. I, 1.4]) and Galois
theory that $E _{\rm sep}$ contains as a subfield an $E$-isomorphic
copy $E _{p} ^{\prime }$ of $E _{p}$, such that ${\cal G}(E _{\rm sep}/(F
_{0}E _{p} ^{\prime }))$ is a Sylow pro-$p$-subgroup of ${\cal G}(E _{\rm
sep}/F _{0})$. Since $F$ is purely inseparable over $F _{0}$, this
implies that ${\cal G}(F _{\rm sep}/(FE _{p} ^{\prime }))$ is a Sylow
pro-$p$-subgroup of ${\cal G}(F _{\rm sep}/F)$, where $F _{\rm sep} = FE
_{\rm sep}$ is the separable closure of $F$ in $\overline E$.
Therefore, $p$ does not divide the degree of any finite extension
of $F$ in $FE _{p} ^{\prime }$. Hence, by (1.2), $D \otimes _{F}
(FE _{p} ^{\prime })$ lies in $d(FE _{p} ^{\prime })$ and has index
ind$(D)$ and exponent exp$(D)$. In addition, it follows that $L
\otimes _{F} (FE _{p} ^{\prime })$ is a cyclic extension of $FE
_{p} ^{\prime }$ which embeds in $D \otimes _{F} (FE _{p} ^{\prime
})$ as an $(FE _{p} ^{\prime })$-subalgebra if and only if $L$
embeds in $D$ over $F$. These observations, statements (1.1)
(ii)-(1.3) and Lemma 8.2 enable one to complete the proof of Lemma
8.3.
\par
\medskip
Now we aim at proving Theorem 8.1 under the hypothesis that $E _{\rm
sep} = E (p)$, for some $p \in {\cal P}(E)$. It is known that ${\cal G} _{E}$ is a
free pro-$p$-group if and only if Br$(E) = \{0\}$ or char$(E) = p$;
this occurs if and only if $H ^{2} ({\cal G} _{E}, \hbox{\Bbb F} _{p}) =
\{0\}$ (cf. [Wa1, Theorem 3.1; Wa2, page 725] or [Se1, Ch. I, 4.2;
Ch. II, 2.2 and 3.1]). Also, it follows from Lemma 4.2 (ii) and
the Albert-Hochschild theorem that if Br$(E) = \{0\}$, then Br$(E
_{1}) = \{0\}$, for every finite extension $E _{1}$ of $E$. For
example, Br$(E) = \{0\}$ when $E$ is perfect and char$(E) = p$
(cf. [A1, Ch. VII, Theorem 22]). Henceforth, we assume that Br$(E)
\neq \{0\}$. Suppose first that $p \neq {\rm char}(E)$. Then Lemma
5.3 implies that $E$ contains a primitive $p$-th root of unity.
Hence, by Lemma 3.8 and Theorem 4.1, $E$ is nonreal and
$p$-quasilocal if and only if ${\cal G} _{E}$ is a $p$-group of Demushkin
type. This completes the proof of Theorem 8.1 (i), so our next
objective is to prove Theorem 8.1 (ii). It follows from [A1, Ch.
XI, Theorem 3] that $E$ is absolutely stable if and only if
$_{p} {\rm Br}(F)$ equals the set $\{[\Delta ]\colon \ \Delta \in
d(F), {\rm ind}(\Delta ) = p\}$, for each finite extension $F/E$.
Note also that central division $F$-algebras of index $p$ are symbol
algebras, since ${\cal G} _{F}$ is a pro-$p$-group. These observations,
combined with (3.2), prove Theorem 8.1 (ii) in
case $p \neq {\rm char}(E)$. In order to complete our proof it
remains to be seen that $E$ is stable, provided that it is almost
perfect, char$(E) = p$ and Br$(E) \neq \{0\}$. This is obtained
from (1.8), [A1, Ch. VII, Theorem 22] and the following lemma.
\par
\medskip
{\bf Lemma 8.4.} {\it Assume that $E$ is a field, such that ${\rm
char}(E) = p > 0$ and $[E\colon E ^{p}] = p$. Let $D \in d(E)$ be
of index $p^{m}$, for some $m \in \hbox{\Bbb N}$. Then ${\rm
exp}(D) = p ^{m}$ and $D$ possesses a maximal subfield that is a
purely inseparable extension of E.}
\par
\medskip
{\it Proof.} Fix an algebraic closure $\overline E$ of $E$ and put
exp$(D) = p ^{\bar m}$, $E _{0} = E$ and $E _{n} = \{\alpha _{n} \in
\overline E\colon \ \alpha _{n} ^{p^{n}} \in E\}$, for every $n \in
\hbox{\Bbb N}$. It follows from (1.8) and the equality $[E\colon E
^{p}] = p$ that $E _{n+1} ^{p} = E _{n}$ and $E _{n}$ is the
unique purely inseparable extension of $E$ in $\bar E$ of degree $p
^{n}$, for each $n \in \hbox{\Bbb N}$. Since, by Albert's theory
of $p$-algebras (cf. [A1, Ch. VII, Theorem 32]), $E _{\bar m}$ is
a splitting field of $D$, this observation and statements (1.1)
(i) and (1.2) (ii) imply that $\bar m = m$ and $E _{m}$ is
$E$-isomorphic to some maximal subfield of $D$, as desired.
\par
\medskip
{\bf Corollary 8.5.} {\it Let $F$ be a quasilocal field, $L/F$ a
finite separable extension and $D \in d(F)$. Then $L$ embeds in
$D$ as an $F$-subalgebra if and only if $[L\colon F]$ divides {\rm
ind}$(D)$; $L$ is a splitting field of $D$ if and only if
$[L\colon F]$ is divisible by {\rm ind}$(D)$.}
\par
\medskip
{\it Proof.} Applying Galois theory, Sylow's theorem and (1.2) as in
the proof of Lemma 8.3, one reduces our considerations to the
special case in which ind$(D)$ is a power of a prime $p$,
$L \subseteq F (p)$ and $L \neq F$. Then our assertion can be
deduced from Theorem 4.1.
\par
\medskip
{\bf Corollary 8.6.} {\it For a quasilocal field $E$, the
following two conditions are equivalent:}
\par
(i) {\it Every finite extension $L$ of $E$ is embeddable as an
$E$-subalgebra in each $\Delta \in d(E)$ of index divisible by
$[L\colon E]$;}
\par
(ii) $E$ {\it has some of the following two properties:}
\par ($\alpha $) $E$ {\it is almost perfect;} ($\beta $) ${\rm
char}(E) = q > 0$ {\it and Br$(E) _{q} = \{0\}$.}
\par
\medskip
{\it Proof.} Suppose first that char$(E) = 0$ or $E$ has property
(ii) ($\beta $). Then char$(E)$ does not divide ind$(\Delta )$,
for any $\Delta \in d(E)$, so it follows from Corollary 8.5 that
condition (i) holds. Henceforth, we assume that char$(E) = q > 0$
and Br$(E) _{q} \neq \{0\}$. The implication (ii) ($\alpha $)$\to
$(i) has essentially been deduced from Corollary 8.5 and Lemma
8.4 in the process of proving [Ch1, Corollary 2.7] (although
formally the result referred to applies to the case where $E$ is
taken as in the concluding assertion of Proposition 2.3). It
remains for us to show here that (i)$\to $(ii) ($\alpha $).
Assuming the opposite, one obtains that there is a purely
inseparable extension $\Phi $ of $E$ such that $[\Phi \colon E] =
q ^{2}$ and $\Phi ^{q} \subseteq E$. At the same time, the
divisibility and nontriviality of Br$(E) _{q}$ guarantees the
existence of an algebra $D \in d(E)$ of exponent $q ^{2}$.
Therefore, by Theorem 3.1 (i), ind$(D) = q ^{2}$. Hence, by [A1,
Ch. VII, Theorem 32], $\Phi $ does not split $D$, which means in
this case that it does not embed in $D$ over $E$. The obtained
contradiction proves that (i)$\to $(ii) ($\alpha $), as required.
\par
\medskip
{\bf Remark 8.7.} Let $F = F _{0} ((X))$ be the formal Laurent power
series field in an indeterminate $X$ over a finite field $F _{0}$ of
characteristic $q$, and let $v$ be the standard $\hbox{\Bbb
Z}$-valued valuation of $F$. It is known that $F$ is noncountable,
hence, $F/F _{0}$ is an extension of infinite transcendency
degree. Fix an infinite set $S _{\infty }$ in $F$ of algebraically
independent elements over $F _{0}$, and for each $n \in
\hbox{\Bbb N}$, denote by $S _{n}$ some subset of $S _{\infty
}$ of cardinality $n$. Also, let $E _{n}$ be the separable
closure of the field $F _{0} (X) (S _{n})$ in $F$, for every $n \in
\hbox{\Bbb N} \cup \{\infty \}$. It is not difficult to see that
the valuation of $E _{n}$ induced by $v$ is Henselian and discrete
with a residue field $F _{0}$. Therefore, by [Ch2, Corollary 2.5],
$E _{n}$ is quasilocal. We show that $E _{n}$ does not possess the
properties of Corollary 8.6 (ii), for any $n \le \infty $. Since
$F _{0}$ has a cyclic extension of degree $q$, $E _{n}$ admits a
nicely semiramified division algebra of index $q$ (see [JW, Sect.
4]), so Br$(E _{n}) _{q} \neq \{0\}$. At the same time, it follows
from the definition of $E _{n}$ and [L1, Ch. X, Propositions 3 and
6] that $[E _{n}\colon E _{n} ^{q}] = q ^{n+1}$, $n \in \hbox{\Bbb
N}$, and $[E _{\infty }\colon E _{\infty } ^{q}] = \infty $.
\par
\medskip
We conclude this Section with examples of quasilocal fields of
very simple type with respect to the structure of the Sylow
subgroups of their absolute Galois groups.
\par
\medskip
{\bf Proposition 8.8.} {\it Let $\overline {\hbox{\Bbb Q}} _{p}$
be an algebraic closure of the field $\hbox{\Bbb Q} _{p}$ of
$p$-adic numbers, $v _{p}$ the unique valuation of $\overline
{\hbox{\Bbb Q}} _{p}$ extending the natural valuation of
$\hbox{\Bbb Q} _{p}$, and $E$ a closed subfield of the completion
$\hbox{\Bbb C} _{p}$ of $\overline {\hbox{\Bbb Q}} _{p}$. Then $E$
is quasilocal.}
\par
\medskip
{\it Proof.} It is well-known that $\hbox{\Bbb C}_{p}$ is
algebraically closed. The assumption on $E$ means that $E$ is
complete with respect to the restriction $v$ of the valuation of
$\hbox{\Bbb C}_{p}$ continuously extending $v_{p}$; this shows in
particular that $\hbox{\Bbb Q}_{p}$ is a subfield of $E$. Observing
that $v(E)$ is a subgroup of $\hbox{\Bbb Q}$, one obtains from the
completeness of $E$ that $v$ is Henselian. Since $\overline
{\hbox{\Bbb Q}} _{p}$ is dense in $\hbox{\Bbb C}_{p}$, this implies
that each finite extension of $E$ in $\hbox{\Bbb C}_{p}$ is
included in an extension of $E$ obtained by adjunction of an
element of $\overline {\hbox{\Bbb Q}} _{p}$ (see the lemma in [L2,
page 380]). Hence, the algebraic closure $\overline E$ of $E$ in
$\hbox{\Bbb C}_{p}$ is equal to $E\overline {\hbox{\Bbb Q}} _{p}$.
In view of Galois theory and the general properties of tensor
products, the obtained result indicates that $\overline E$ is
$E$-isomorphic to $E \otimes _{E _{0}} \overline {\hbox{\Bbb Q}}
_{p}$, where $E _{0} = E \cap \overline {\hbox{\Bbb Q}} _{p}$.
Moreover, it becomes clear that every $E$-automorphism of
$\overline E$ is determined by its action on $\overline
{\hbox{\Bbb Q}} _{p}$, and also, that ${\cal G} _{E}$ can be identified
with ${\cal G} _{E _{0}}$. As $\hbox{\Bbb Q} _{p}$ is quasilocal,
Proposition 8.8 follows now directly from Theorem 8.1 and Lemma
8.2.
\par
\medskip
{\bf Proposition 8.9.} {\it With notation being as in Theorem
8.1, let $E$ be a nonreal field and $G_{p}$ a pro-$p$-group of rank
$n(p) \le 2$, for every $p \in \Pi (E)$. Then $E$ is quasilocal.}
\par
\medskip
{\it Proof.} Lemma 8.3 allows one to consider only the special case
in which $E _{\rm sep} = E (p)$, for some $p \in {\cal P}(E)$. It follows
from [J, Proposition 4.4.8], Lemma 4.2 (ii) and the
Albert-Hochschild theorem that if $p = {\rm char}(E)$, then Br$(U)
= \{0\}$, for every algebraic extension $U/E$. Assuming that $p
\neq {\rm char}(E)$ whence $E$ contains a primitive $p$-th root of
unity, one obtains from [Wa2, Lemma 7] that ${\cal G} _{E}$ is a free
pro-$p$-group or a Demushkin group (see also [EnV, Theorem 4.7],
for the case of $p = 2$). Now our assertion can be deduced from
Lemma 3.8 and Theorem 4.1.
\par
\medskip
Note finally that the conditions of Proposition 8.9 hold, if $E$
is a quasilocal field of characteristic $q \ge 0$ with some of the
following two properties: (i) ${\cal G} _{E}$ is torsion-free with abelian
Sylow pro-$p$-subgroups, for all $p \in \Pi (E)$ (apply Lemma 8.3
and [Ch3, Lemma 3.2]); (ii) $q \not\in \Pi (E)$ and $E$ has a
Henselian valuation $v$ such that $v(E)$ is $p$-indivisible, for
every $p \in \Pi (E)$ (cf. [Ch1, (1.2) (ii) and Remark 2.2] and
[Ch2]).
\par
\vskip1.75truecm \centerline{\bf Acknowledgements}
\par
\medskip
I would like to thank the referee for the extremely careful reading
of the paper and for his/her suggestions, especially, the
improvements in the presentation of Sections 1, 5 and Remark 8.7,
the correction of (6.2), the observed inaccuracies in Sections 6 and
7, and the much shorter proofs of Proposition 2.1 and the former
part of Lemma 4.2 (i). \vskip1cm \centerline{ REFERENCES}
\vglue15pt\baselineskip12.8pt
\def\num#1{\smallskip\item{\hbox to\parindent{\enskip [#1]\hfill}}}
\parindent=1.38cm
\par
\medskip
\num{A1} A.A. {\pc ALBERT}, {\sl Structure of Algebras.} Amer. Math.
Soc. Coll. Publ., 24, Amer. Math. Soc., XII, New York, 1939.
\par
\num{A2} A.A. {\pc ALBERT}, {\sl Modern Higher Algebra.} Chicago
Univ. Press, XIV, Chicago, Ill., 1937.
\par
\num{Am1} S.A. {\pc AMITSUR}, {\sl Generic splitting fields of
central simple algebras.} Ann. Math. 62 (1955), 8-43.
\par
\num{Am2} S.A. {\pc AMITSUR}, {\sl On central division algebras.}
Isr. J. Math. 12 (1972), 12-22.
\par
\num{Am3} S.A. {\pc AMITSUR}, {\sl Generic splitting fields.} Proc.
Brauer Groups in Ring Theory and Algebraic Geometry (Antwerp 1981),
Lecture Notes in Math. 917 (1981), 1-24.
\par
\num{Ar} M. {\pc ARTIN}, {\sl Two-dimensional orders of finite
representation type.} Manuscripta math. 58 (1987), 445-471.
\par
\num{B} E. {\pc BRUSSEL}, {\sl  Division algebra subfields
introduced by an indeterminate.} J. Algebra 188 (1997), 216-255.
\par
\num{CF} J.W.S. {\pc CASSELS}, A. {\pc FR$\ddot o$HLICH} (Eds.),
{\sl Algebraic Number Theory.} Proc. Instruct. Conf., organized by
LMS (a NATO Adv. Study Inst.) with the support of IMU, held at the
Univ. of Sussex, Brighton, 01.9-17.9., 1965, Academic Press,
London-New York, 1967.
\par
\num{Ch1} I.D. {\pc CHIPCHAKOV}, {\sl Henselian valued stable
fields.} J. Algebra 208 (1998), 344-369.
\par
\num{Ch2} I.D. {\pc CHIPCHAKOV}, {\sl Henselian valued quasilocal
fields with totally indivisible value groups.} Comm. Algebra 27
(1999), No 7, 3093-3108.
\par
\num{Ch3} I.D. {\pc CHIPCHAKOV}, {\sl On the Galois cohomological
dimensions of Henselian valued stable fields.} Comm. Algebra 30
(2002), 1549-1574.
\par
\num{Ch4} I.D. {\pc CHIPCHAKOV}, {\sl Central division algebras of
$p$-primary dimensions and the $p$-component of the Brauer group of
a $p$-quasilocal field.} C.R. Acad. Sci. Bulg. 55 (2002), 55-60.
\par
\num{Ch5} I.D. {\pc CHIPCHAKOV}, {\sl On the residue fields of
Henselian valued stable fields.} Preprint, Inst. Math. Bulg. Acad.
Sci., February 1997, No 1, 21 p.
\par
\num{Ch6} I.D. {\pc CHIPCHAKOV}, {\sl One-dimensional abstract local
class field theory.} Preprint, arXiv:math/0506515v4 [math.RA].
\par
\num{Ch7} I.D. {\pc CHIPCHAKOV}, {\sl On the Brauer groups of
quasilocal nonreal fields and the norm groups of their finite Galois
extensions.} Preprint, arXiv:0707.4245v3 [math.RA].
\par
\num{Ch8} I.D. {\pc CHIPCHAKOV}, {\sl Primarily quasilocal fields,
divisible torsion Galois modules and Brauer groups of stable
fields.} Preprint.
\par
\num{D1} S.P. {\pc DEMUSHKIN}, {\sl On the maximal $p$-extension of
a local field.} Izv. Akad. Nauk SSSR, Math. Ser., 25 (1961), 329-346
(in Russian).
\par
\num{D2} S.P. {\pc DEMUSHKIN}, {\sl On $2$-extensions of a local
field.} Sibirsk. Mat. Zh. 4 (1963), 951-955 (Russian: English
transl. in Amer. Math. Soc. Transl., Ser. 2, 50 (1966), 178-182).
\par
\num{Dr1} P. {\pc DRAXL}, {\sl Skew Fields.} London Math. Soc.,
Lecture Note Series, 81. Cambridge etc., Cambridge Univ. Press,
1983.
\par
\num{Dr2} P. {\pc DRAXL}, {\sl Ostrowski's theorem for Henselian
valued skew fields.} J. Reine Angew. Math. 354 (1984), 213-218.
\par
\num{Ef} I. {\pc EFRAT}, {\sl On fields with finite Brauer groups.}
Pac. J. Math. 177 (1997), 33-46.
\par
\num{E} O. {\pc ENDLER}, {\sl Valuation Theory.} Springer-Verlag,
Berlin-Heidelberg-New York, 1972.
\par
\num{EnV} A.J. {\pc ENGLER}, T.M. {\pc VISWANATHAN}, {\sl Digging
holes in algebraic closures a la Artin I.} Math. Ann. 265 (1983),
263-271.
\par
\num{Er} Yu.L. {\pc ERSHOV}, {\sl Co-Henselian extensions and
Henselizations of division algebras.} Algebra i Logika 27 (1988),
649-658 (Russian: English transl. in Algebra and Logic 27 (1988),
401-407).
\par
\num{FSS} B. {\pc FEIN}, D. {\pc SALTMAN}, M. {\pc SCHACHER}, {\sl
Heights of cyclic field extensions.} Bull. Soc. Math. Belg. 40
(1988), Ser. A, 213-223.
\par
\num{FS} B. {\pc FEIN}, M. {\pc SCHACHER}, {\sl Relative Brauer
groups, I.} J. Reine Angew. Math. 321 (1981), 179-194.
\par
\num{FSa} T.J. {\pc FORD}, D. {\pc SALTMAN}, {\sl Division algebras
over Henselian surfaces.} Israel Math. Conf. Proc., 1989, No 1,
320-336.
\par
\num{F} L. {\pc FUCHS}, {\sl Infinite Abelian Groups.} Academic
Press, New York-London, 1970.
\par
\num{JW} B. {\pc JACOB}, A. {\pc WADSWORTH}, {\sl Division algebras
over Henselian fields.} J. Algebra 128 (1990), 126-179.
\par
\num{J} N. {\pc JACOBSON}, {\sl Finite-Dimensional Division Algebras
over Fields.} Springer-Verlag, Berlin, 1996.
\par
\num{Jo} A.J. de {\pc JONG}, {\sl The period-index problem in the
theory of algebraic surfaces.} Duke Math. J. 123 (2004), 71-94.
\par
\num{K} G. {\pc KARPILOVSKY}, {\sl Topics in Field Theory.}
North-Holland Math. Stud. 155, Amsterdam, 1989.
\par
\num{Ko} H. {\pc KOCH}, {\sl Galoissche Theorie der
$p$-Erweiterungen.} Springer-Verlag, New York-Heidelberg-Berlin,
1970.
\par
\num{Lab1} J.P. {\pc LABUTE}, {\sl Demushkin groups of rank $\aleph
_{0}$.} Bull. Soc. Math. France 94 (1966), 211-244.
\par
\num{Lab2} J.P. {\pc LABUTE}, {\sl Classification of Demushkin
groups.} Canad. J. Math. 19 (1967), 106-132.
\par
\num{La} T.Y. {\pc LAM}, {\sl Orderings, valuations and quadratic
forms.} Conf. Board Math. Sci. Regional Conf. Ser. Math. No 52,
Amer. Math. Soc., Providence, RI, 1983.
\par
\num{L1} S. {\pc LANG}, {\sl Algebra.} Addison-Wesley, Reading, MA,
1965.
\par
\num{L2} S. {\pc LANG}, {\sl On quasi-algebraic closure.} Ann. Math.
55 (1952), 373-390.
\par
\num{LvdD} A. {\pc LUBOTZKY}, L. {\pc VAN DEN DRIES}, {\sl Subgroups
of free profinite groups and large subfields of Q.} Isr. J. Math. 39
(1981), 25-45.
\par
\num{M} A.S. {\pc MERKURJEV}, {\sl Brauer groups of fields.} Comm.
Algebra 11 (1983), No 22, 2611-2624.
\par
\num{MS} A.S. {\pc MERKURJEV}, A.A. {\pc SUSLIN}, {\sl
$K$-cohomology of Severi-Brauer varieties and norm residue
homomorphisms.} Izv. Akad. Nauk SSSR 46 (1982), 1011-1046 (Russian:
English transl. in Math. USSR Izv. 21 (1983), 307-340).
\par
\num{MW1} J. {\pc MINA$\breve c$}, R. {\pc WARE}, {\sl Demushkin
groups of rank $\aleph _{0}$ as absolute Galois groups.} Manuscripta
math. 73 (1991), 411-421.
\par
\num{MW2} J. {\pc MINA$\breve c$}, R. {\pc WARE}, {\sl
Pro-$2$-Demushkin groups as Galois groups of maximal $2$-extensions
of fields.} Math. Ann. 292 (1992), 337-353.
\par
\num{P} R. {\pc PIERCE}, {\sl Associative Algebras.} Graduate Texts
in Mathematics, 88, Springer-Verlag, New York-Heidelberg-Berlin,
1982.
\par
\num{Pl} V.P. {\pc PLATONOV}, {\sl The Tannaka-Artin problem and
reduced $K$-theory.} Izv. Akad. Nauk SSSR 40 (1976), 227-261
(Russian: English transl. in Math. USSR Izv. 10 (1976), 211-243
(1977)).
\par
\num{Re} I. {\pc REINER}, {\sl Maximal Orders.} London Math. Soc.
Monographs, v. 5, London-New York-San Francisco: Academic Press, a
subsidiary of Harcourt Brace Jovanovich, Publishers, 1975.
\par
\num{R} P. {\pc RIBENBOIM}, {\sl Equivalent forms of Hensel's
lemma.} Exposition. Math. 3 (1985), 3-24.
\par
\num{Roq1} P. {\pc ROQUETTE}, {\sl On the Galois cohomology of the
projective linear group and its applications to the construction of
generic splitting fields of algebras.} Math. Ann. 150 (1963),
411-439.
\par
\num{Roq2} P. {\pc ROQUETTE}, {\sl Isomorphism of generic splitting
fields of simple algebras.} J. Reine Angew. Math. 214 (1964),
207-226.
\par
\num{Se1} J.-P. {\pc SERRE}, {\sl Cohomologie Galoisienne.} Lecture
Notes in Math. 5, Springer-Verlag, Berlin-Heidelberg-New York, 1965.
\par
\num{Se2} J.-P. {\pc SERRE}, {\sl Structure de certains
pro-$p$-groupes (apr$\grave e$s Demushkin).} Semin. Bourbaki 15
1962/1963 (1964), No 252, 11 p.
\par
\num{Se3} J.-P. {\pc SERRE}, {\sl Local Fields.} Graduate Texts in
Mathematics, 67, Springer-Verlag, New York-Heidelberg-Berlin, 1979.
\par
\num{T} J.-P. {\pc TIGNOL}, {\sl On the corestriction of central
simple algebras.} Math. Z. 194 (1987), 267-274.
\par
\num{TY} I.L. {\pc TOMCHIN}, V.I. {\pc YANCHEVSKIJ}, {\sl On defects
of valued division algebras.} Algebra i Analiz 3 (1991), No 3,
147-164 (Russian: English transl. in St. Petersburg Math. J. 3
(1992), No 3, 631-646).
\par
\num{W} A.R. {\pc WADSWORTH}, {\sl Extending valuations to finite
dimensional division algebras.} Proc. Amer. Math. Soc. 98 (1986),
20-22.
\par
\num{Wa1} R. {\pc WARE}, {\sl Quadratic forms and profinite
$2$-groups.} J. Algebra 5 (1979), 227-237.
\par
\num{Wa2} R. {\pc WARE}, {\sl Galois groups of maximal
$p$-extensions.} Trans. Amer. Math. Soc. 333 (1992), No 2, 721-728.
\par
\num{Wh} G. WHAPLES, {\sl Algebraic extensions of arbitrary fields.}
Duke Math. J. 24 (1957), 201-204. \vskip1cm
\def\pc#1{\eightrm#1\sixrm}
\hfill\vtop{\eightrm\hbox to 5cm{\hfill Ivan {\pc CHIPCHAKOV}\hfill}
 \hbox to 5cm{\hfill Institute of Mathematics and Informatics\hfill}\vskip-2pt
 \hbox to 5cm{\hfill Bulgarian Academy of Sciences\hfill}
\hbox to 5cm{\hfill Acad. G. Bonchev Str., bl. 8\hfill} \hbox to
5cm{\hfill 1113 {\pc SOFIA,} Bulgaria\hfill}}
\end
\par
\bye